\RequirePackage{fix-cm}

 \documentclass[showkeys,preprint,nofootinbib]{revtex4-1}          
\usepackage[utf8]{inputenc} 
\usepackage{lmodern,textcomp,graphicx}

\usepackage{amsmath,amssymb}

\makeatletter
\@namedef{u8:\detokenize{α}}{\ensuremath{\alpha}}
\@namedef{u8:\detokenize{β}}{\ensuremath{\beta}}
\@namedef{u8:\detokenize{Γ}}{\ensuremath{\Gamma}}
\@namedef{u8:\detokenize{γ}}{\ensuremath{\gamma}}
\@namedef{u8:\detokenize{Δ}}{\ensuremath{\Delta}}
\@namedef{u8:\detokenize{δ}}{\ensuremath{\delta}}
\@namedef{u8:\detokenize{ε}}{\ensuremath{\epsilon}}
\@namedef{u8:\detokenize{ϵ}}{\ensuremath{\varepsilon}}
\@namedef{u8:\detokenize{ζ}}{\ensuremath{\zeta}}
\@namedef{u8:\detokenize{η}}{\ensuremath{\eta}}
\@namedef{u8:\detokenize{Θ}}{\ensuremath{\Theta}}
\@namedef{u8:\detokenize{θ}}{\ensuremath{\theta}}
\@namedef{u8:\detokenize{ϑ}}{\ensuremath{\vartheta}}
\@namedef{u8:\detokenize{ι}}{\ensuremath{\iota}}
\@namedef{u8:\detokenize{κ}}{\ensuremath{\kappa}}
\@namedef{u8:\detokenize{Λ}}{\ensuremath{\Lambda}}
\@namedef{u8:\detokenize{λ}}{\ensuremath{\lambda}}
\@namedef{u8:\detokenize{μ}}{\ensuremath{\mu}} 
\@namedef{u8:\detokenize{µ}}{\ensuremath{\mu}} 
\@namedef{u8:\detokenize{ν}}{\ensuremath{\nu}}
\@namedef{u8:\detokenize{ξ}}{\ensuremath{\xi}}
\@namedef{u8:\detokenize{Π}}{\ensuremath{\Pi}}
\@namedef{u8:\detokenize{∏}}{\ensuremath{\prod}}
\@namedef{u8:\detokenize{π}}{\ensuremath{\pi}}
\@namedef{u8:\detokenize{ϖ}}{\ensuremath{\varpi}}
\@namedef{u8:\detokenize{ρ}}{\ensuremath{\rho}}
\@namedef{u8:\detokenize{ϱ}}{\ensuremath{\varrho}}
\@namedef{u8:\detokenize{Σ}}{\ensuremath{\Sigma}}
\@namedef{u8:\detokenize{σ}}{\ensuremath{\sigma}}
\@namedef{u8:\detokenize{ς}}{\ensuremath{\varsigma}}
\@namedef{u8:\detokenize{τ}}{\ensuremath{\tau}}
\@namedef{u8:\detokenize{υ}}{\ensuremath{\upsilon}}
\@namedef{u8:\detokenize{Υ}}{\ensuremath{\Upsilon}}
\@namedef{u8:\detokenize{Φ}}{\ensuremath{\Phi}}
\@namedef{u8:\detokenize{φ}}{\ensuremath{\varphi}}
\@namedef{u8:\detokenize{ϕ}}{\ensuremath{\phi}}
\@namedef{u8:\detokenize{χ}}{\ensuremath{\chi}}
\@namedef{u8:\detokenize{Ψ}}{\ensuremath{\Psi}}
\@namedef{u8:\detokenize{ψ}}{\ensuremath{\psi}}
\@namedef{u8:\detokenize{Ω}}{\ensuremath{\Omega}}
\@namedef{u8:\detokenize{ω}}{\ensuremath{\omega}}

\@namedef{u8:\detokenize{∑}}{\sum}
\@namedef{u8:\detokenize{∇}}{\ensuremath{\nabla}}
\@namedef{u8:\detokenize{⦅}}{\ensuremath{\left(}}
\@namedef{u8:\detokenize{⦆}}{\ensuremath{\right)}}
\@namedef{u8:\detokenize{×}}{\ensuremath{\times}}

\@namedef{u8:\detokenize{ｄ}}{\ensuremath{\mathrm{d}}}
\@namedef{u8:\detokenize{⟨}}{\ensuremath{\langle}}
\@namedef{u8:\detokenize{⟩}}{\ensuremath{\rangle}}
\@namedef{u8:\detokenize{⟪}}{\ensuremath{\llangle}}
\@namedef{u8:\detokenize{⟫}}{\ensuremath{\rrangle}}
\@namedef{u8:\detokenize{∂}}{\ensuremath{\partial}}
\@namedef{u8:\detokenize{·}}{\ensuremath{\cdot}}
\@namedef{u8:\detokenize{÷}}{\frac}
\@namedef{u8:\detokenize{⋀}}{\ensuremath{\wedge}}
\@namedef{u8:\detokenize{ℝ}}{\ensuremath{\mathbb{R}}}
\newsavebox{\@brx}
\newcommand{\llangle}[1][]{\savebox{\@brx}{\(\m@th{#1\langle}\)}%
  \mathopen{\copy\@brx\mkern2mu\kern-0.9\wd\@brx\usebox{\@brx}}}
\newcommand{\rrangle}[1][]{\savebox{\@brx}{\(\m@th{#1\rangle}\)}%
  \mathclose{\copy\@brx\mkern2mu\kern-0.9\wd\@brx\usebox{\@brx}}}
\makeatother

\newcommand\vt{\mathbf}
\newcommand\tsr{\mathbb}
\renewcommand\u{\vt u}

\newcommand{\td}[2] {\cfrac{\textrm{d} #1}{\textrm{d} #2}\,}
\newcommand{\pd}[2] {\cfrac{\partial #1}{\partial #2}\,}
\newcommand\di\displaystyle
\usepackage{ifthen}
\newcommand{\ppd}[3]{
   \ifthenelse{\equal{#2}{#3}}
   {\cfrac{\partial^2 #1}{\partial {#2}^2}\,}         
   {\cfrac{\partial^2 #1}{\partial #2\partial #3}\,}  
}

\newcommand\p{\vt p}
\newcommand\q{\vt q}
\renewcommand\v{\vt v}

\newcommand\ov\overline
\newcommand\e{\textrm{e}}
\usepackage{afterpage}

\newcommand\tp{^{\mathsf{T}}}

\newcounter{subsubsubsection}[subsubsection]

\usepackage{amsthm}
\newtheorem{theorem}{Theorem}
\newtheorem{corollary}[theorem]{Corollary}
\newtheorem{lemma}[theorem]{Lemma}

\begin{document}


\title{Some robust integrators for large time dynamics}

\author{Dina Razafindralandy}
\email{dina.razafindralandy@univ-lr.fr}
\affiliation{LaSIE -- University of La Rochelle, 
Av. Michel Cr\'epeau, 17042 La Rochelle Cedex 1, France}

\author{Vladimir Salnikov} 
\email{vladimir.salnikov@univ-lr.fr}
\affiliation{LaSIE  -- CNRS \& University of La Rochelle,
Av. Michel Cr\'epeau, 17042 La Rochelle Cedex 1, France}

\author{Aziz Hamdouni}
\email{aziz.hamdouni@univ-lr.fr}
\affiliation{LaSIE -- University of La Rochelle, 
Av. Michel Cr\'epeau, 17042 La Rochelle Cedex 1, France}

\author{Ahmad Deeb}
\email{ahmad.deeb@univ-lr.fr}
\affiliation{LaSIE -- University of La Rochelle, 
Av. Michel Cr\'epeau, 17042 La Rochelle Cedex 1, France}

\begin{abstract} 
This article reviews some integrators particularly suitable for the numerical resolution of differential equations on a large time interval. Symplectic integrators are presented. Their stability on exponentially large time is shown through numerical examples. Next, Dirac integrators for constrained systems are exposed. An application on chaotic dynamics is presented. Lastly, for systems having no exploitable geometric structure, the Borel-Laplace integrator is presented. 
Numerical experiments on Hamiltonian and non-Hamiltonian systems are carried out, as well as on a partial differential equation.


\end{abstract}

\maketitle



\section{Introduction}

In many domains of mechanics, simulations over a large time interval are crucial. This is, for instance, the case in molecular dynamics, in weather forecast or in astronomy. While many time integrators are available in literature, only few of them are suitable for large time simulations. Indeed, many numerical fail to correctly predict the expected physical phenomena such as energy preservation,  as the simulation time grows.

For equations having an underlying geometric structure (Hamiltonian systems, variational problems, Lie symmetry group, Dirac structure, \dots), geometric integrators appear to be very robust for large time simulation. These integrators mimic the geometric structure of the equation at the discrete scale. 

The aim of this paper is to present some time integrators which are suitable for large time simulations. We consider not only equations having a geometric structure but also more general equations. We first present symplectic integrators for Hamiltonian systems. We show in section \ref{section:symplectic} their ability in preserving the Hamiltonian function and some other integrals of motion. Applications will be on a periodic Toda lattice and on $n$-body problems. To simplify, the presentation is done in canonical coordinates. 

In section \ref{section:dirac}, we show how to fit a constrained problem into a Dirac structure. We then detail how to construct a geometric integrator respecting the Dirac structure. The presentation will be simplified, and the (although very interesting) theoretical geometry is skipped. References will be given for more in-depth understanding. A numerical experiment, showing the good long time behavior of Dirac integrators, will be carried out.

In section \ref{section:borel}, we present the Borel-Padé-Laplace integrator (BPL). BPL is a general-purpose time integrator, based on a time series decomposition of the solution, followed by a resummation to enlarge the validity of the series and then reducing the numerical cost on a large time simulation. Finally, the long time behaviour will be investigated through numerical experiments on Hamiltonian and non-Hamiltonian system, as well as on a partial differential equation. Numerical cost will be examined when relevant.

\section{Symplectic integrators\label{section:symplectic}}

We first make some reminder on Hamiltonian systems and their flows in canonical coordinates. Some examples of symplectic integrators are given afterwards and numerical experiments are presented.

\subsection{Hamiltonian system}

A Hamiltonian system in $ℝ^d×ℝ^d$ is a system of differential equations which can be written as follows:
\begin{equation}
	\begin{cases}
		\td {\q}t=\ \pd H{\p},\\\\
		\td {\p}t=-\pd H{\q},
	\end{cases}
	\label{hamiltonian0}
\end{equation}
the Hamiltonian $H$ being a function of time and of the unknown vectors $\q=(q_1,..,q_d)\tp$ and $\p=(p_1,\cdots,p_d)\tp$.
Equations (\ref{hamiltonian0}) can be written in a more compact way as follows:
\begin{equation}
	\td{\u}t=\tsr J\nabla H
	\label{hamiltonian_J}
\end{equation}
where $\u=(\q,\p)\tp$, $\nabla H=\pd H{\u}$ and $\tsr J$ is the skew-symmetric matrix
\[ \tsr J=\begin{pmatrix}
	0&\tsr I_{d}\\-\tsr I_{d}&0
\end{pmatrix},
\]
$\tsr I_d$ being the identity matrix of  $ℝ^{d}$. 

The flow of the Hamiltonian system (\ref{hamiltonian_J}) at time $t$ is the function $Φ_t$ which, to an initial condition $\u^0$ associates the solution $\u(t)$ of the system. More precisely, $Φ_t$ is defined as:
\begin{equation}Φ_t:\ 
	\begin{array}{rcl}
		ℝ^d×ℝ^d&\longrightarrow&ℝ^d×ℝ^d\\\\
		\u^0=(\q^0,\p^0)\tp&\longmapsto&\u(t)=(\q(t),\p(t))\tp.
	\end{array}
	\label{flow}
\end{equation}
The property that $\tsr J^{-1}=\tsr J\tp=-\tsr J$ leads to the symplecticity property of $Φ_t$:
\begin{equation}
	(\nabla Φ_t)\tp\ \tsr J\ (\nabla Φ_t)=\tsr J.\label{symplecticity}
\end{equation}
Note that $\tsr J$ can be seen as an area form, in the following sense. If $\vt v$ and $\vt w$ are two vectors of $ℝ^d×ℝ^d$, with components 
\[\vt v=(v_{q_1},\cdots,v_{q_d},v_{p_1},\cdots,v_{p_d})\tp,\quad\vt w=(w_{q_1},\cdots,w_{q_d},w_{p_1},\cdots,w_{p_d})\tp\] then 
\[\vt v\tp\ \tsr J\ \vt w=\sum_{i=1}^d⦅v_{q_i}w_{p_i}-v_{p_i}w_{q_i}⦆.\]
In words, $\vt v\tp\ \tsr J\ \vt w$ is the sum of the areas formed by $\vt v$ and $\vt w$ in the planes $(q_i,p_i)$. The symplecticity property (\ref{symplecticity}) then means that the flow of a Hamiltonian system is area preserving.

In the sequel, $H$ is assumed autonomous in time. It can then be shown that $H$ is preserved along trajectories.

\subsection{Flow of a numerical scheme}

Consider a numerical scheme which computes an approximation $\u^n$ of the solution $\u(t^n)$ of equation (\ref{hamiltonian_J}) at time $t^n$. The flow of this scheme is defined as 
\begin{equation}φ_{t^n}:\quad\quad
		\u^0=(\q^0,\p^0)\tp\quad\longmapsto\quad\u^n=(\q^n,\q^n)\tp
	\label{flow_discrete}
\end{equation}
For a one-step integrator, with a time step $Δt=t^{n+1}-t^n$ (which may depend on $n$), it is more convenient to work with the one-step flow
\begin{equation}φ_{Δt}^n: \quad\quad
		\u^n\quad\longmapsto\quad\u^{n+1}
	\label{flow_discrete1}
\end{equation} 
As an example, consider the explicit Euler integration scheme
\[ 
\begin{cases}
	\q^{n+1}=\q^n+{Δt}\,\pd H{\p}(\q^n,\p^n)
	\\\\
	\p^{n+1}=\p^n-{Δt}\,\pd H{\q}(\q^n,\p^n).
\end{cases}
\]
The one-step flow of this scheme is
\begin{equation}
	φ_{Δt}^n(\u^n)=\u^n+{Δt}\,\tsr J∇H(\u^n).
\end{equation}

The one-step flow of a fourth order Runge-Kutta scheme is
\begin{equation}
	φ_{Δt}^n(\u^n)=\u^n+Δt÷{f_1+2f_2+2f_3+f_4}6
	\label{RK4_flow}
\end{equation}
where 
\[f_1=\tsr J∇H(\u^n), \quad \quad \quad\quad \quad \quad  f_2=\tsr J∇H(\u^n+\tfrac{Δt}2 f_1), \]
\[f_3=\tsr J∇H(\u^n+\tfrac{Δt}2 f_2), \quad \quad \quad f_4=\tsr J∇H(\u^n+Δt\ f_3).\]

\subsection{Some symplectic integrators}

A time integrator is called symplectic if its flow is symplectic, meaning that
\begin{equation}
	(∇ φ_{t^n})\tp\ \tsr J\ (∇ φ_{t^n})=\tsr J.
\end{equation}
Geometrically, a symplectic integrator is then a time scheme which preserves the area form. For a one-step scheme, this property is equivalent to
\begin{equation}
	(∇ φ_{Δt}^n)\tp\ \tsr J\ (∇ φ_{Δt}^n)=\tsr J\label{symplecticity_one}
\end{equation}
at each iteration $n$. 

When $d=1$, it can easily be shown that, for an explicit Euler scheme,
\begin{equation}
	(∇ φ_{Δt}^n)\tp\ \tsr J\ (∇ φ_{Δt}^n)=⦅1+\Delta t^2⦅\ppd Hq q \ppd Hpp-\ppd Hqp⦆⦆\tsr J,
\end{equation}
meaning that the explicit Euler scheme is not symplectic. Neither the implicit Euler scheme is symplectic.
 By contrast, by mixing the explicit and the implicit Euler scheme, we get a symplectic scheme, called symplectic Euler scheme, defined as follows
\begin{equation}
	\begin{cases}
		\q^{n+1}=\q^n+Δt\pd H{\p}(\q^n,\p^{n+1}),
		\\\\
		\p^{n+1}=\p^n-Δt\pd H{\q}(\q^n,\p^{n+1}).
	\end{cases}
	\label{symplectic_euler}
\end{equation}
Note that in (\ref{symplectic_euler}), one can take $(\q^{n+1},\p^n)$ in the right-hand-side instead of $(\q^n,\p^{n+1})$. This leads to the other symplectic Euler scheme. 

Both symplectic Euler schemes are first order. A way to get a second order scheme is to compose two symplectic Euler schemes with time steps $Δt/2$. One then obtains the St\"ormer-Verlet schemes \cite{hairer03}. Another way is to take the mid-points in the right-hand side of an Euler scheme \cite{hairer06}.

Symplectic Runge-Kutta schemes of higher order can be built as follows. An $s$-stage Runge-Kutta integrator of equation (\ref{hamiltonian0}) is defined as \cite{hairer_87,sanz92}:
\begin{equation}
	\begin{array}{l}\displaystyle
		\vt U_i=\u^n+Δt\ \sum_{j=1}^sα_{ij}\ \tsr J∇H(\vt U_j), \quad i=1,\cdots,s, 
		\\\\\displaystyle
		\u^{n+1}=\u^n+Δt\ \sum_{i=1}^sβ_i\ \tsr J∇H(\vt U_i),
	\end{array}
	\label{rk}
\end{equation}
for some real numbers $α_{ij}$, $β_i$, $i,j=1,\cdots,s$. Scheme (\ref{rk}) is symplectic if and only if the coefficients verify the relation (\cite{lasagni88,sanz88}):
\begin{equation}
	β_iβ_j=β_iα_{ij}+β_jα_{ji},\quad\quad i,j=1,\cdots,s.
\end{equation}
An example of symplectic Runge-Kutta scheme is the 4th order, 3-stage scheme defined by the coefficients
\begin{equation}
	α=\begin{pmatrix}
		\di÷b2&&&0&&&0 \\\\b&&&\di÷12-b&&&0\\\\b&&&1-2b&&&\di÷b2
	\end{pmatrix},\quad \quad
	β=\begin{pmatrix}
		b&&&1-2b&&&b
	\end{pmatrix},
	\label{rk4sym}
\end{equation}
where
\[
	b=÷{2+2^{1/2}+2^{-1/3}}3.
\]
Many other variants of symplectic Runge-Kutta methods can be found in the literature (see for example \cite{feng10}). 

Symplectic integrators do not preserve exactly the Hamiltonian in general. However, the symplecticity condition seems to be strong enough that, experimentally, symplectic integrators exhibit a good behaviour toward the preservation property. In fact, we have the following error estimation on the Hamiltonian \cite{benettin1994,hairer97}:
\begin{equation}
	|H(p^n,q^n)-H(p^0,q^0)|=O(Δt^r)\quad\quad \textrm{ for }nΔt\leq\e^{÷{γ}{2Δt}}\label{hamiltonian_bounded}
\end{equation}
for some constant $γ>0$, $r$ being the order of the symplectic scheme. This relation states that the error is bounded over an exponentially long discrete time. Moreover, it was shown in \cite{cooper87} that symplectic Runge-Kutta methods preserve exactly quadratic invariants.

In the next subsection, some interesting numerical properties of symplectic schemes are highlighted through some model problems.

\subsection{Numerical experiments}


\subsubsection{Periodic Toda lattice\label{section:toda}}

The evolution of a periodic Toda lattice with $d$ particles can be described with the Hamiltonian function
\[ H=\sum_{k=1}^d⦅÷12p_k^2+\e^{q_k-q_{k+1}}⦆ \]
where $q_k$ is the (one-dimensional) position of the $k-$th particle, $q_{d+1}=q_1$ and $p_k$ is its momentum. A periodic Toda lattice is completely integrable and  it is known that the eigenvalues of the following matrix $L$ are first integrals of the system:
\begin{equation}
	L=\begin{pmatrix}
		a_1&b_1&0&0&0&\cdots&0&b_d\\
		b_1&a_2&b_2&0&0&\cdots&0&0\\
		0&b_2&a_3&b_3&0&\cdots&0&0\\
		\vdots\\
		0&0&0&\cdots&0&b_{d-2}&a_{d-1}&b_{d-1}\\
		b_d&0&0&\cdots&0&0&b_{d-1}&a_d\\
	\end{pmatrix}
	\label{l}
\end{equation}
where
\[ a_k=-÷12p_k,\quad\quad b_k=÷12 \e^{÷12(q_k-q_{k+1})}.\]

In the numerical test, we consider $d=3$ particles, positionned initially at $q_1=0$, $q_2=2$ and $q_3=3$. The initial momenta are $p_1=0.5$, $p_2=-1.5$ and $p_3=1$. 

\begin{figure}[ht]
	\includegraphics[width=7cm]{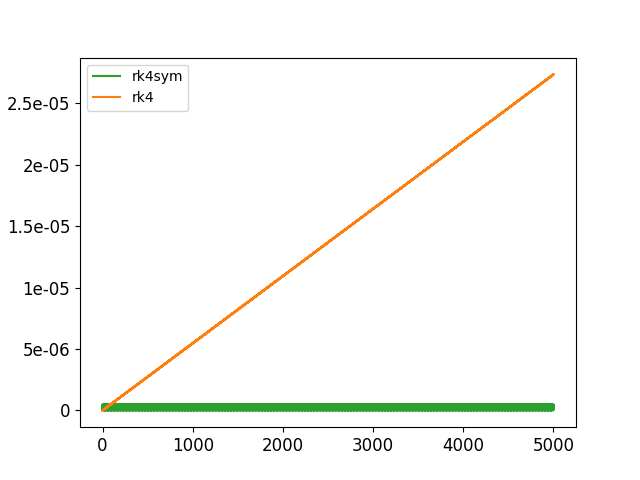}
	\caption{Toda Lattice. Relative error $÷{|H(\u^n)-H(\u^0)|}{|H(\u^0)|}$ with $Δt=10^{-2}$}
	\label{fig_toda01_errH}
\end{figure}

\begin{figure}[ht]
	\includegraphics[width=60mm]{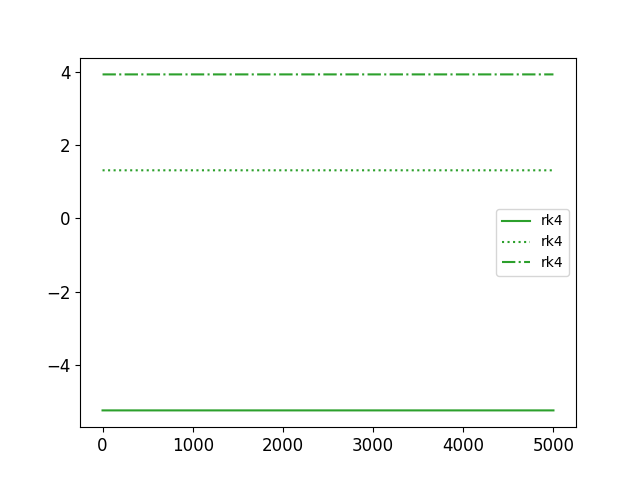}
	\includegraphics[width=60mm]{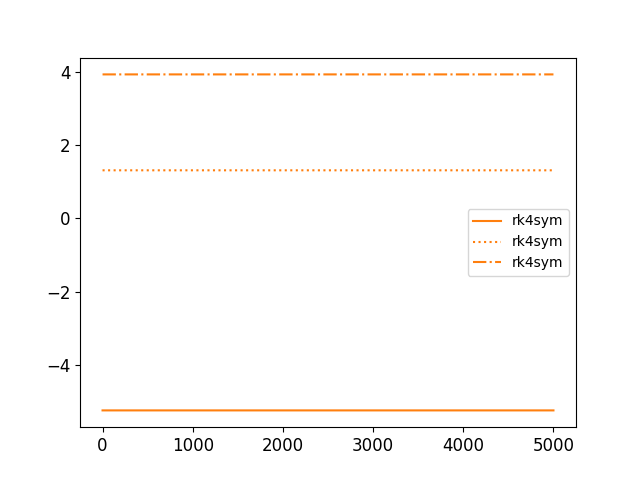}
	\caption{Toda Lattice. Evolution of the eigenvalues of $L$ with $Δt=10^{-2}$. Left: RK4. Right: RK4sym}
	\label{fig_toda01_eigen}
\end{figure}

First, we choose a time step $Δ=10^{-2}$. The Hamiltonian equation is solved with the classical 4-th order Runge-Kutta scheme (RK4) and the symplectic version (RK4sym) defined by (\ref{rk4sym}) up to $t=5000$. The relative error on the Hamiltonian is presented on figure \ref{fig_toda01_errH}. As can be seen, the RK4 error oscillates and increases globally. It remains acceptable for $t<5000$ since it does not exceed $2.735·10^{-5}$. The RK4sym error also oscillates but is much closer to zero. Its highest value is about $4.625·10^{-7}$, that is an order of $10^{-2}$ bellow RK4 error at $t=5000$. Both schemes globally preserve the three eigenvalues of the matrix $L$, as can be observed on figure~\ref{fig_toda01_eigen}.

\begin{figure}[ht]
	\includegraphics[width=6cm]{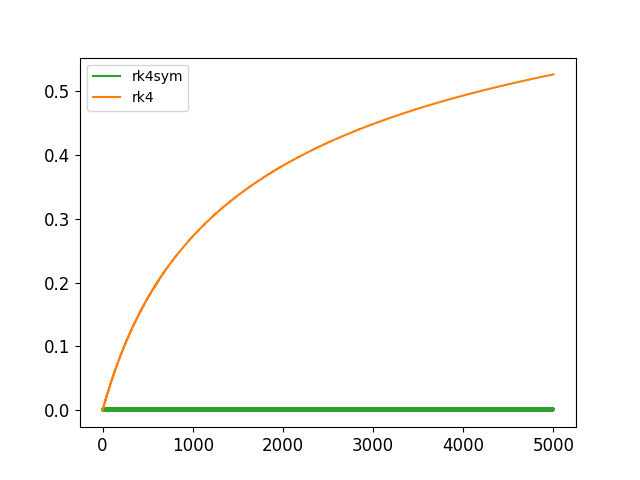}
	\includegraphics[width=6cm]{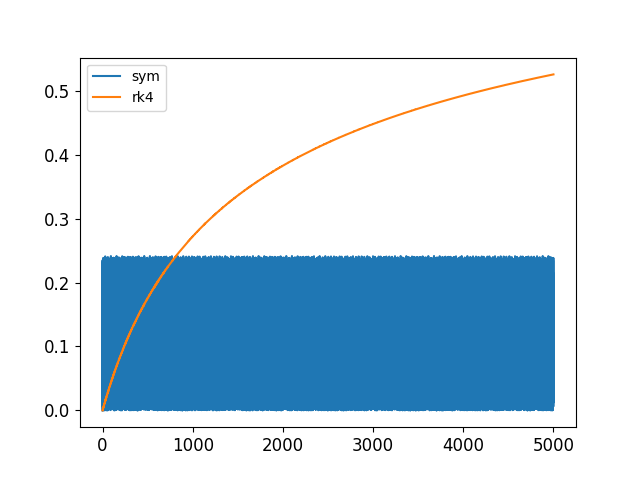}
	\caption{Toda Lattice. Relative error on $H$ with $Δt=10^{-1}$. Left: RK4 and RK4sym. Right: RK4 and Symplectic Euler}
	\label{fig_toda1_errH}

\end{figure}

Next, $Δt$ is set to $10^{-1}$. With this time step, the error of the classical Runge-Kutta scheme increases quickly from the first iterations, as can be observed in figure \ref{fig_toda1_errH}. It reaches 50 percent around $t=4.2·10^{3}$. As for it, the RK4sym error oscillates around $2.26·10^{-3}$ but does not present any increasing global tendency. The error given by the symplectic Euler scheme (\ref{symplectic_euler}) is also plotted on figure \ref{fig_toda1_errH}, left. It oscillates around $0.83·Δt$ and does not exceed $2.42Δt$. So, for $t$ greater than 820, even the first order symplectic Euler scheme produces an error smaller than the 4-th order non-symplectic Runge-Kutta scheme on the Hamiltonian.

\begin{figure}[ht]
	\includegraphics[width=6cm]{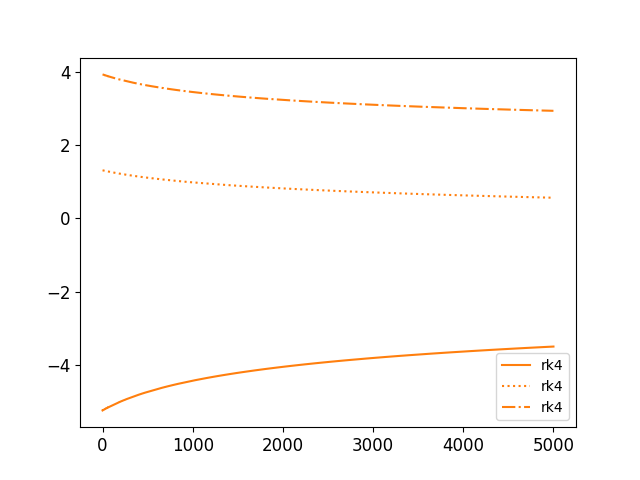}
	\includegraphics[width=6cm]{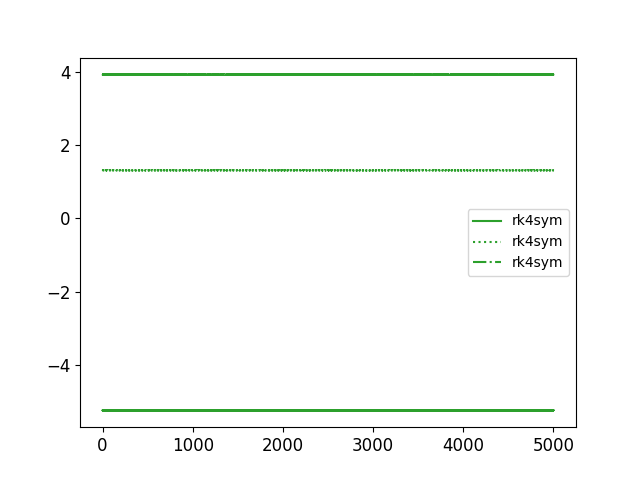}
	\caption{Toda Lattice. Evolution of the eigenvalues of $L$ with $Δt=10^{-1}$. Left: RK4. Right: RK4sym}
	\label{fig_toda1_eigen}
\end{figure}

\begin{figure}[ht]
	\includegraphics[width=6cm]{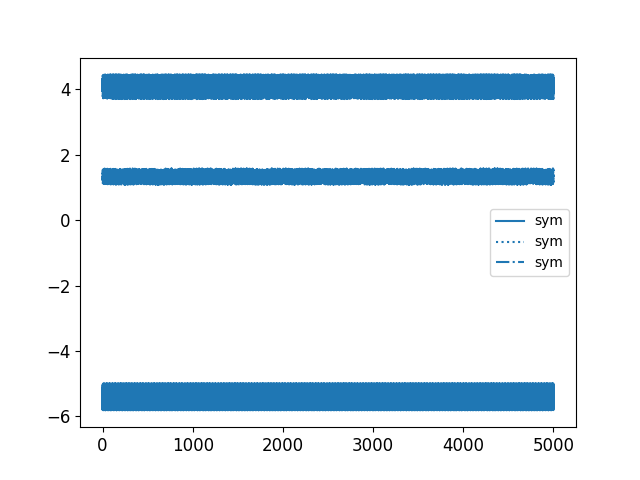}
	\caption{Toda Lattice. Evolution of the eigenvalues of $L$ with the symplectic Euler scheme and $Δt=10^{-1}$.}
	\label{fig_toda1_eigen_sym}
\end{figure}

The evolution of the eigenvalues of the matrix $L$ is presented on figure \ref{fig_toda1_eigen}. As can be seen, the eigenvalues are not preserved by the classical Runge-Kutta scheme. For example, the computed smallest eigenvalue at $t=5000$ is about $-3.50$ whereas its initial value is $-5.24$. With the symplectic Runge-Kutta scheme, the error on the smallest eigenvalue oscillates around $4.85·10^{-3}$ but does not present an increasing tendency.
With the first order symplectic Euler scheme, the oscillations are much more pronounced, as can be observed on figure \ref{fig_toda1_eigen_sym}, but as with RK4sym, there is no increasing trend. It becomes smaller than the RK4 error when the simulation time increases.

It is clear from these experiments that symplectic schemes are particularly stable for large time simulations, where the user wishes a time step as large as possible to reduce the computation cost. In some situation, even a symplectic scheme with a smaller order gives better results over a long time than a classical integrator.

\subsection{$n$-body problem}

As a second example, consider the system of $d$ bodies subjected to mutual gravitational forces. The evolution of the system is described by the Hamiltonian function
\[ H=\sum_{k=1}^{d}÷12÷{\|\p_k\|^2}{m_k}-\sum_{k=1}^{d}\sum_{l=k+1}^d÷{Gm_km_l}{\|\q_l-\q_k\|} .\]
In this expression, $\q_k$ and $\p_k$ are the position vector and the momentum of the $k-$th body, $m_k$ is its mass and $G$ is the gravitation constant.

\begin{figure}[ht]
	\centering
	\includegraphics[width=6cm]{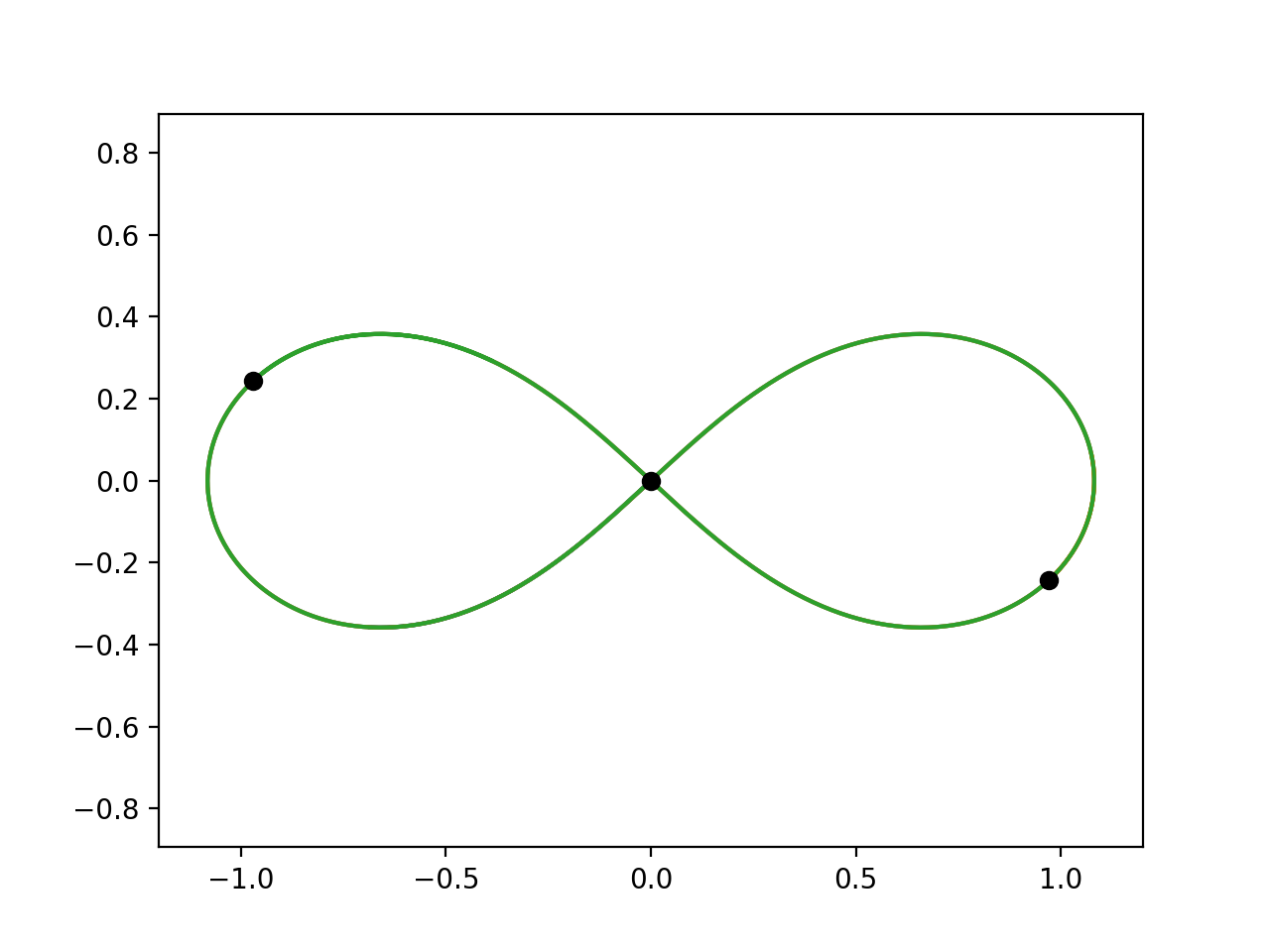}
	\includegraphics[width=6cm]{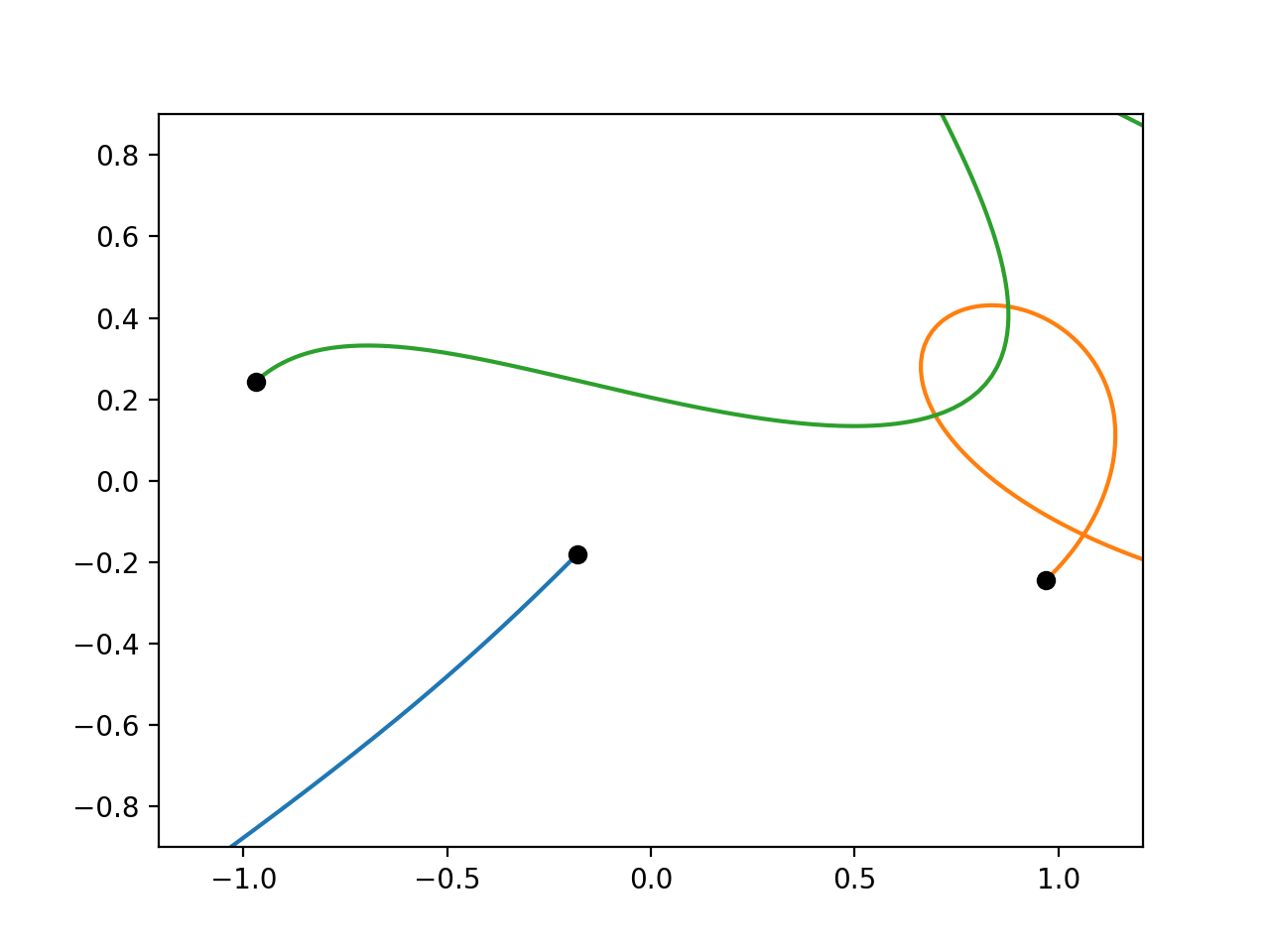}
	\caption{Three-body problem. Figure-eight orbit (left) and perturbed (right) initial positions and trajectories}
	\label{fig_3body_eight.png}
\end{figure}

\begin{figure}[ht]

	\includegraphics[width=6cm]{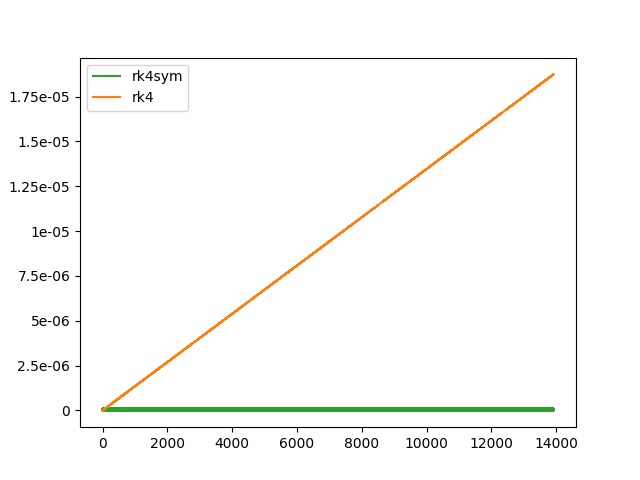}
	\includegraphics[width=6cm]{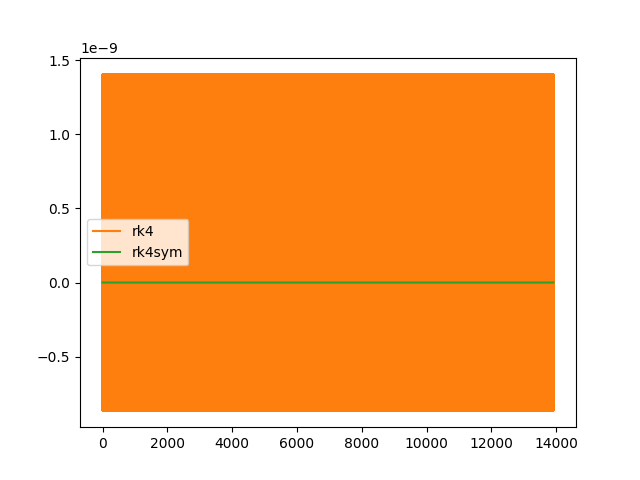}
	\caption{Eight-orbit configuration. Error on the Hamiltonian (left) and on the angular momentum}
	\label{fig_3body_l}

\end{figure}

For the simulation, we take $d=3$. We consider the initial configuration corresponding to the choregraphic figure-eight in \cite{chenciner00}. The solution is periodic, with period $T\simeq 6.32591398$. The common trajectory and the initial position are presented on figure \ref{fig_3body_eight.png}, left. The simulation is run up to $t=2200T$, with a time step $\Delta t=0.02T$. In these configurations, the classical Runge-Kutta scheme provides a fairly good result but the symplectic scheme is much more accurate regarding the preservation properties. The relative errors on the Hamiltonian and the error on the angular momentum are plotted on figure \ref{fig_3body_l}. The error on $H$ is about $1.87·10^{-5}$ at the final time with RK4 and $4.375·10^{-10}$ with RK4sym.

\begin{figure}[ht]
	\includegraphics[width=6cm]{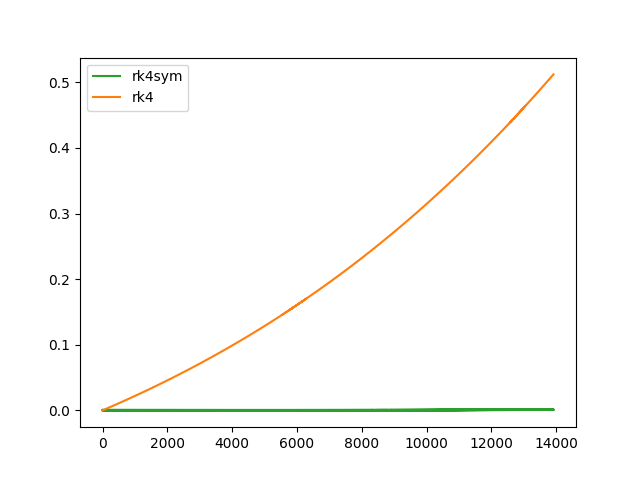}
	\includegraphics[width=6cm]{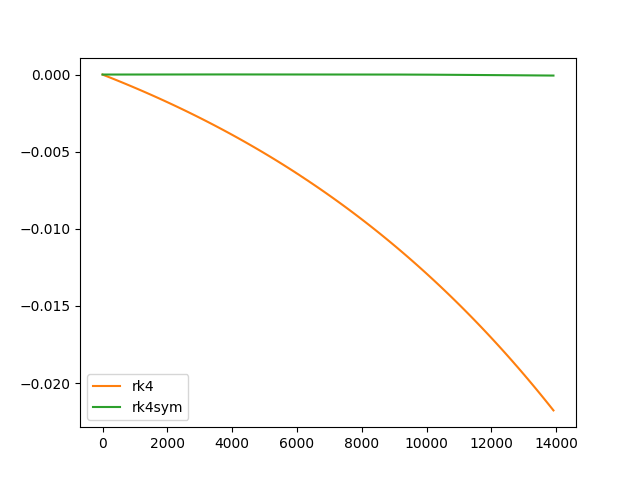}
	\caption{Perturbed configuration. Error on the Hamiltonian (left) and on the angular momentum}
	\label{fig_3body_pl_errH_pl}
\end{figure}

In a second simulation, the initial position and momentum of the body in the middle in figure \ref{fig_3body_eight.png} are changed into their values at $t=T/80$ in the figure-eight solution. The initial configuration of the bodies at the two ends are kept as in the previous simulation. In this case, the figure-eight is broken. A part of the trajectorie of each body and their initial positions are presented on figure \ref{fig_3body_eight.png}, right. The evolution of the error on the Hamiltonian and on the angular momentum is plotted on figure \ref{fig_3body_pl_errH_pl}. As can be noticed, the error on the Hamiltonian increases quickly with the classical RK4. It reaches 50 percent at $t\simeq2166T\simeq13702$. At this time value, the error on the Hamiltonian with the symplectic RK4 is about $7.75·10^{-4}$. It reaches 1 percent much later, around $t\simeq6315T\simeq39950$.

These numerical experiments show again that symplectic schemes are more robust than classical ones for long time dynamics simulation. They have a good behaviour regarding the preservation of the Hamiltonian and some other first integrals. Similar results have been obtained in a previous work (\cite{amses18}) on an harmonic oscillator, on Kepler's problem and on vortex dynamics.

Obviously, not all mechanical systems fit into the Hamiltonian formalism, hence there are other geometric constructions that are worth being considered in the context of structure-preserving integrators.

\section{Dirac integrators\label{section:dirac}}

In this section we give an overview of the so-called Dirac structures, and describe a 
class of mechanical systems where those appear naturally, namely systems with constraints. Originally, Dirac structures appear in the work of T.~Courant (\cite{courant}). The initial motivation was coming from mechanics. As is known, for mechanical systems one can choose between Lagrangian and Hamiltonian formalisms, both being equivalent in finite dimension. The rough idea behind Dirac 
structures is to consider both formalisms simultaneously, i.e. working with velocities \emph{and}
momenta, however not forgetting that those are dependent variables. 
Geometrically, this means that instead of choosing between the tangent and cotangent bundles $TM$ or $T^*M$ for the phase space, we consider their direct sum $E = TM \oplus T^*M$ and a subbundle of it, subject to some compatibility conditions.  Somehow, the original work did not have direct applications to mechanics, since the geometry of the problem turned out to be rather intricate, and gave rise to 
a lot of development in \emph{higher structures} and in \emph{theoretical physics}. 
In the last decade, however, it was revived with the introduction of so-called port--Hamiltonian  (\cite{port-ham}) and implicit Lagrangian systems (\cite{YoMa1, YoMa2}).

\subsection{Geometric construction}

Not to overload the presentation here with geometric details, let us talk about spaces instead of bundles. To recover the geometric picture, a motivated reader is invited to read the original paper \cite{courant} of Courant or the overview of relevant results in \cite{dirac-int}.
The object of study will be $\mathbb{R}^{2d}$ (in the same notations as in the beginning
of the previous section), and some natural construction around it.
 We are going to view it as 
$\mathbb{R}^{2d} = \mathbb{R}^d \times V^*$, that is the trivial bundle 
over $\mathbb{R}^d$ with a fiber being $V^*$ -- the dual of some $d$-dimensional vector space $V$.  Morally, $V$ is the space of velocities $\v$ at each point $\q$ and $V^*$ corresponds to the space of momenta $\p$. In coordinates: 
 \[\q = (q_1, \dots, q_d)\tp \in \mathbb{R}^d, \quad
\v =  (v_1, \dots, v_d)^T \in V, \quad
\p = (p_1, \dots, p_d)^T \in V^*. \]

In this setting, imposing constraints on the system  means defining some 
restriction on couples $(\q, \v)$, i.e. not all the points $\q$ are permitted, and at each point 
$\q$, $\v$ is not arbitrary, but belongs to a subspace of $V$. Under some 
regularity conditions, one can say, that at each velocity $\v$ belongs to a set $\Delta_{\q}\subset V$ which is the kernel of a set of 
linear forms $\alpha^a$. And since everything depends on the point $\q$, globally, 
the permitted vector fields $\v(\q)$ live in the kernel $\Delta\subset ℝ^d×V$ of $m$ differential $1$-forms 
$\alpha^a(\q), a = 1, \dots, m$.

To transfer this construction from $\mathbb{R}^d \times V$ to $\mathbb{R}^d \times V^*$,  one  needs to consider 
double vector bundles (\cite{tul}). In our simplified setting this means that the space of interest is ${\cal V} = \mathbb{R}^{4d}$, where each component has some geometric 
interpretation. Namely, we consider ${\cal V}$ as the tangent to $\mathbb{R}^d \times V^*$. Naturally, $\mathbb{R}^d \times V$ is embedded in ${\cal V}$ (recall that $V$ is  tangent to $\mathbb{R}^d$). The constraint set is then a subset $\tilde \Delta \subset {\cal V}$, and the 
differential forms $\alpha^a(\q)$ generate its annihilator $\Delta_0$ that naturally belongs to 
${\cal V}^*$. 
Note that, since $\mathbb{R}^d \times V^*$ is a symplectic space. 
It is equipped with a bilinear antisymmetric non-degenerate closed form $\Omega$ (this form $\Omega$ is the generalisation of the matrix $\tsr J$ of section \ref{section:symplectic} in non-canonical coordinates). One can then construct 
a symplectic mapping $\Omega^\flat \colon {\cal V} \to {\cal V}^*$. 

We can now define the \emph{Dirac structure}\footnote{The correct terminology is
\emph{almost} Dirac structure. For details see \cite{dirac-int}.} associated to the system with constraints:
\begin{equation} \label{DDeltaQ}
\mathbb{D}_{\Delta}  = \{ (w, \beta)  \in {\cal V} \times {\cal V}^* \;| \;   
  w \in \tilde \Delta, \;
\beta - \Omega^\flat w \in \Delta_0    \}.  \nonumber 
\end{equation} 

To define the dynamics of the system, we introduce two more objects.
First, the Lagrangian, which as usual is  a mapping 
$L \colon \mathbb{R}^d\times V \to \mathbb{R}$ induces its differential 
which is a mapping $d L$ from $\mathbb{R}^d\times V$ to its cotangent.
And again by post-composing it with a symplectomorphism of the appropriate double 
bundles, one constructs the  \emph{Dirac differential} which locally reads
\[{\cal D} L\colon 
	\begin{array}{rcl}
		\mathbb{R}^d\times V &\longrightarrow& {\cal V}^*
		\\\\
		(\q, \v) &\longmapsto& \displaystyle\left(\q,\ \frac{\partial L}{\partial \v},\ -\frac{\partial L}{\partial \q},\ \v\right)
	\end{array}
\]
Second, the evolution of  the system will be described by a so-called 
 \emph{partial vector field} $X$, i.e. a mapping 
 \[X \; \colon \;  \Delta \oplus Leg(\Delta)\ \subset\ \left(\mathbb{R}^d\times V\right) \oplus \left(\mathbb{R}^d\times V^*\right) \quad\longrightarrow\quad {\cal V},\] 
where $Leg(\Delta)$ is the image of $\Delta$ by the Legendre transform. 
$X$ should be viewed as a vector field on $\mathbb{R}^d\times V^*$, 
with the momenta parametrized by the Legendre transform of the velocities compatible with the constraints.

With the above notations, the \emph{implicit Lagrangian system} 
 is a triple $(L, \Delta, X)$, such that $(X, {\cal D}L) \in \mathbb{D}_{\Delta}$.
In local coordinates, this means: 
\begin{equation} \label{impl_lagr}
	  \begin{cases}
		  \td {\q}t \in \Delta, \quad\quad\displaystyle
		  \p = \frac{\partial L}{\partial \v}, \\\\
		  \td {\q}t = \v,   \quad\quad\displaystyle
		  \td{\p}t - \frac{\partial L}{\partial \q} \in \Delta_0.
	  \end{cases}
\end{equation}
One understands easily the mechanical interpretation of the first three equations 
above. The forth one can be rewritten as 
\[
\td{\p}t - \frac{\partial L}{\partial q} =  \sum_{a=1}^m \lambda_a \alpha^a,
\]
and one recognizes immediately the Lagrange multipliers.

\subsection{Discretization}
It is important to note that the previous section is not just a ``fancy'' way of recovering the well-known theory: every step  of the construction admits a discrete 
analog. We briefly present the recipe  of this discretization and, again, refer the interested reader to \cite{dirac-int} for details and examples.

 The system is characterized by the following continuous data:
the Lagrangian function $L(q, \dot q)$ and the set of constraint 1-forms 
$\alpha^a, a = 1, \dots, m$.  
The discrete version $L_d$ of the Lagrangian at time $t^n$ is 
\[
  L_d = Δt\ L(\q^n, \v^n).
\]
And the constraints are rewritten as
\[
 < \alpha_d^a , \v^n> = 0,\quad a = 1, \dots, m.
\]
where $< \alpha_d^a , \v^n>=\alpha_d^a(\v^n)$, and $α_d^a$ is a discrete version of $α^a$.
In the equations above, $\q^n$ is the value of $\q$ and $\v^n$ is an approximation of the velocity $\v$, both at time $t^n$.
 
 To construct the numerical method out of these data, one applies the following procedure:
\begin{eqnarray}
  \p^{n+1} &=& \frac{1}{Δt}\, \frac{\partial L_d}{\partial \v^n} \label{line1} \\\nonumber\\
  \p^n - \frac{1}{Δt}\, \frac{\partial L_d}{\partial \v^n} +  \frac{\partial L_d}{\partial \q^n} &=& \sum_{a=1}^m \lambda_a 
  \frac{\partial < \alpha_d^a, \v^n>}{\partial \v^n} \label{line2} \\\nonumber\\
  < \alpha_d^a, \v^n> &=& 0, \quad a= 1, \dots, m. \label{line3}
\end{eqnarray}
The variables appearing explicitly are the values of $\p$ at the $n$-th and the $(n+1)$-st steps, and $\v^n$ should be some approximation of the velocity, thus 
bringing $\q^n$ to the system.
Here, we consider two natural options: 
\begin{itemize}
	\item $\v^n := \cfrac{\q^{n+1} - \q^n}{Δt}$,\quad we label it Dirac-1, and

	\item $\v^n := \cfrac{\q^{n+1} - \q^{n-1}}{2Δt}$, \quad labelled Dirac-2.
\end{itemize}
Dirac-1 permits to recover the method from \cite{leok} whereas we introduced Dirac-2 in \cite{dirac-int}. In both cases, we obtain $2d + m$ equations: $d$ from each of the lines 
 (\ref{line1}) and (\ref{line2}) of the equations above, and $m$ from the constraints
 (\ref{line3}). At the $n$-th step the unknowns are $\q^{n+1}$, $\p^{n+1}$ and $\boldsymbol{\lambda}=(\lambda_1,\dots,\lambda_m)$, so the system we obtain is complete. It is linear in $\boldsymbol{\lambda}$ and $\p^{n+1}$,
  and when the constraints are holonomic, in $\q^{n+1}$ as well.  

It is important to note that, in some sense, this construction generalizes the previous section.
If one considers the system without constraints but still applies the procedure, 
(\ref{line3}) becomes obsolete, and in (\ref{line2}) the right-hand-sides vanish, 
so one obtains a numerical method for the dynamics of a Lagrangian system governed by $L$. By a straightforward computation, one checks that for a natural mechanical system with a potential $U$, i.e. when 
 $L = \frac{1}{2}m\v^2 + U(\q)$, Dirac-1 is symplectic.
 And it is also meaningful to consider a symplectic version of Dirac-2 (we do not detail it here since we would need to explain  what is symplecticity  for a multistep method).

\subsection{Example: chaos for double pendulum}

\begin{figure}[ht]
	\includegraphics[width=27mm]{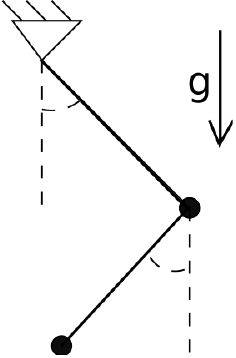}
	\caption{double pendulum}
	\label{fig:2pen}
\end{figure}

We will apply the Dirac integrators constructed in the previous subsection to a scholar problem of a planar double pendulum in a gravity field: 
a system of 2 mass points attached to rigid inextensible weightless rods (see figure \ref{fig:2pen}).
Although, it looks like a classical well-studied problem, it is a bit challenging for simulations. In the absence of the gravity field, this is a textbook  example of an integrable system (energy and angular momentum are conserved, thus the Liouville--Arnold theorem can be applied). With the gravity, the ``folkloric knowledge'' 
says that it is chaotic, although as far as we know there is no rigorous proof of the fact. For the integrabilty, there is a semi-numerical  proof of absence of an additional first integral in \cite{2pen}, using the computation of the monodromy group by a method presented in \cite{monodromy}. Anyway, the apparent chaoticity of the system
results in its sensitivity to parameters and initial data in the numerical simulation.

From the point of view of the previous subsection, this is a typical example 
of a system with constraints: the distance $\ell_1$ from the first mass point to the origin and 
the distance  $\ell_2$ between the two mass points are fixed. 
The system admits a parametrization in terms of angles, but we will pretend not 
to  know it, to test the method.

We thus consider a mechanical system of two mass points given by the Lagrangian
\[
L(\q_1,\q_2, \dot {\q}_1, \dot {\q}_2) = \frac1{2}m_1 \|\dot {\q}_1\|^2 + \frac1{2}m_2 \|\dot {\q}_2\|^2  - m_1gq_{1,y} - m_2gq_{2,y},
\]
subject to the constraints 
\[
   \varphi^1 \equiv \|\q_1\|^2 - \ell_1^2 = 0, \quad\quad    \varphi^2 \equiv \|\q_2-\q_1\|^2 - \ell_2^2 = 0.
   \]
To recover the framework of the numerical method given by (\ref{line1} - \ref{line3}), we 
take $\alpha^a \equiv ｄ\varphi^a, a = 1,2$.

\begin{figure}[htp]
	\includegraphics[width=10cm]{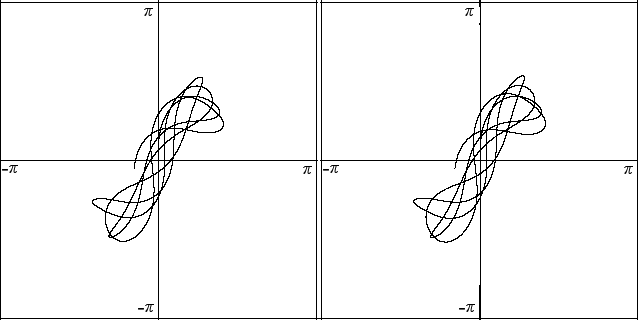}
	\includegraphics[width=10cm]{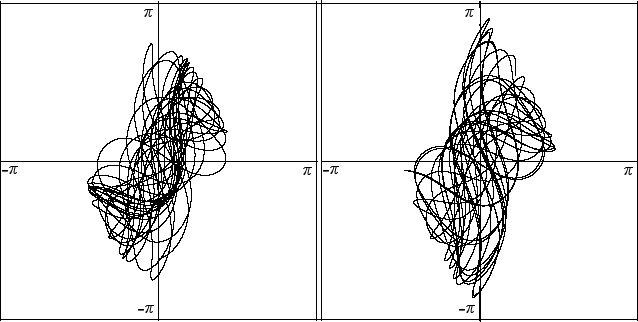}
	\includegraphics[width=10cm]{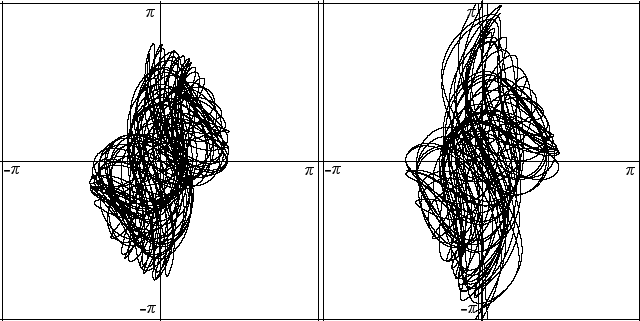}

	\caption{Double pendulum: Comparison of dynamics. Left: Dirac-2, right: explicit Euler. \newline
	Region swept by the trajectory until $T = 10, \; 50 , \; 100$, respectively, from top to bottom. }
	\label{fig:euler-vs-dirac}
\end{figure}

The typical result of simulations is shown on figure \ref{fig:euler-vs-dirac}.
 Dirac-2  and explicit Euler methods are compared. For visualization (but not for computation),
 we use the angle representation of the double pendulum (see figure \ref{fig:2pen}).
Both algorithms start with the same initial data, and the same 
timestep $\Delta t = 0.0001$. They are in good agreement in the beginning as can be seen on the two top-graphics of figure \ref{fig:euler-vs-dirac}. But already at time $T = 50$, the difference is visible (graphics in the middle). And towards 
$T=100$ the difference becomes dramatic: for the Euler method the pendulum is making 
full turns instead of oscillation. And this is clearly a computation artifact, since decreasing the timestep one gets rid of the discrepancy and recovers the left picture for both methods.
Note also that Dirac structure based method preserves the constraints much better than the 
Euler one: the error is  $2.2\cdot 10^{-6}$ compared to  $0.06$ respectively. 

A similar effect is observed 
for other methods: trapezium, and even Runge--Kutta, which is of higher order.
Moreover, there is another non-negligible convenience of the Dirac structure based methods: 
the Lagrange multipliers are treated like other dynamical variables, there is no need to resolve ``by hand'' the equation related to constraints (\ref{line3}).

\vspace{\baselineskip}

In many areas of mechanics, systems are often described by an underlying geometric structure. As observed in the two previous sections, making use of these structures leads to more robust numerical schemes. In the last section, we propose an integrator which is suitable to general systems where no geometric structure is exploitable for numerical simulations.

\section{Borel-Laplace integrator\label{section:borel}}

Consider an ordinary differential or a semi-discretized partial differential equation:
\begin{equation}
	\td{u}t=F(u(t),t). \label{eq}
\end{equation}
We look for a time series solution:
\begin{equation}
	\breve{u}(t)=\sum_{n=0}^{+\infty}u_nt^n. \label{series}
\end{equation}
Generally, the terms $u_n$ of the series can be computed via a recurrence relation of the form
\begin{equation}
	(n+1)u^{n+1}=F^n(u_0,\dots,u_n) .
	\label{recurrence}
\end{equation}
In (\ref{recurrence}), $F^n$ is the $n$-th Taylor expansion of the function $F$.
A Borel summation is applied to series (\ref{series}). This summation is essential if the convergence radius of the series is zero and the series is summable \cite{borel99,ramis91}. If the convergence radius is not zero, the Borel summation enlarges the domain of validity of the series. 

The Borel sum of $\breve{u}(t)$ is 
\begin{equation}
	\mathcal{S}\breve{u}(t)=[\mathcal{L}\circ\mathcal{P}\circ \mathcal{B}]\breve{u}(t)
	\label{borel}
\end{equation}
where $\mathcal{B}$ is the Borel transform, $\mathcal{P}$ is a prolongation along a semi-line in $\mathbb{C}$ linking 0 to infinity (we will take the real positive semi-line), and $\mathcal{L}$ is the Laplace transform along this semi-line. The theory of Borel summation can be found, for example, in \cite{borel99,ramis12a,ramis12b,ramis91}. Some other works on BPL, as a time integrator, can be found in \cite{jcp13,esaim14}. In this section, we present very briefly the Borel-Padé-Laplace algorithm, integrated into a numerical scheme.

The Borel-Padé-Laplace summation integrator (BPL) consists in the following steps:
\begin{itemize}
	\item Given an initial condition $u(t_0)=u_0$, compute a truncated series solution via recurrence (\ref{recurrence}):
		$\displaystyle
			\breve{u}^N(t)=\sum_{n=0}^{N}u_nt^n 
		$
		.
	\item Compute its Borel transform:
		$\displaystyle\mathcal{B}\breve{u}^N(ξ)=\sum_{n=0}^{N-1}÷{u_{n+1}}{n!}\,ξ^n.$
	\item Transform $\mathcal{B}\breve{\u}^N(ξ)$ into a rational fraction function via a Padé approximation: $\displaystyle P^N(ξ)=÷{a_0+a_1t+\dots a_{N_{num}}t^{N_{num}}}{b_0+b_1t+\dots b_{N_{den}}t^{N_{den}}}$

		 The Padé approximation materializes the prolongation in the Borel summation procedure. 
	\item Apply a Laplace transformation (at $1/t$) on $P(ξ)$ to obtain a numerical Borel sum
		$\displaystyle\mathcal{S}\breve{u}^N(t)=\int_0^{+\infty}P^N(\xi)\e^{-\xi/t}ｄ\xi.$
		
		Numerically, the integral is replaced by a Gauss-Laguerre quadrature.

	\item Take $\mathcal{S}\breve{u}^N(t)$ as an approximate solution $u(t)$ of (\ref{eq}) within the integral $[t_0,t_1]$ where the residue of the equation is smaller than a parameter $ε_{res}$. 
	\item Restart the algorithm with $u_0=u(t_1)$ as initial condition to obtain an approximate solution for larger values of $t$.
\end{itemize}
At each iteration, $t_1-t_0$ is considered as the (adaptative) time-step of the scheme. The average time-step will be used for comparisons in numerical experiments. Note that at each time, the approximate solution has an analytical representation as a Laplace integral. A continued fraction representation can also be used \cite{dcds16}.

An advantage of BPL is that it is totally explicit, in contrast with symplectic integrators in general. Moreover, changing the order of the scheme is as easy as setting $N$ to a different value. Note also that the resummation procedure can be done componentwise, enabling an easy parallezation on multi-core computers. However, no such optimization has been done in the present article.

In the following subsection, an attempt to study the symplecticity property of BPL is presented.

\subsection{High-order symplecticity}

The numerical flow of BPL can be defined as
\[ u_0\quad\mapsto\quad \mathcal{S}\breve{u}^N(t).\]
Currently, no symplecticity result has yet been found on this scheme. Instead, it can be shown that a scheme based on the truncated series $\breve{u}$ without the resummation procedure is symplectic at order $N$, if the equation is symplectic. 

The flow of the scheme based on the time series $\breve{u}$ is
\[ φ_{t,\breve u}:\quad u_0\quad\mapsto\quad \breve{u}(t)=\sum_{n=0}^{+\infty}t^nu_n.\]
\begin{lemma}
	The flow of the scheme based on the time series, applied to the Hamiltonian equation (\ref{hamiltonian_J}) is
	\begin{equation}
		φ_{t,\breve u}=∑_{n=0}^{+\infty}÷{t^n}{n!}\tsr JD^n∇H \quad\quad\text{where}\quad D^n=\cfrac{ｄ^{n-1}}{ｄt^{n-1}}
		\label{phiseries}
	\end{equation}
	agreeing that $D^0=1$.
\end{lemma}
This can be straightforwardly deduced by injecting the time series $\breve{u}$ in (\ref{hamiltonian_J}) and identifying the coefficients of each $t^n$. Next, if the series is convergent then, inside the convergence disc, $\breve{u}$ is the exact solution. In this case, $\breve{φ}_t$ is symplectic. We reformulate this statement in the following theorem.
\begin{theorem}
	If the series is convergent then
	\begin{equation}
		(∇φ_{t,\breve u})\tp\ \tsr J\ φ_{t,\breve u}=\tsr J. 
		\label{symp}
	\end{equation}
\end{theorem}
\begin{corollary}
	If the series is convergent then, for any $n\geq1$,
	\begin{equation}
		∑_{k=0}^n ÷1{k!}÷1{(n-k)!}(\tsr JD^k∇H)\tp\ \tsr J\ (\tsr JD^{n-k}∇H)=0.
		\label{crl}
	\end{equation}
\end{corollary}
This corollary is obtained by injecting the series development (\ref{phiseries}) into (\ref{symp}) and identifying the coefficients of $t^n$ for $n\geq1$. For $n=0$, we simply have
\begin{equation}
	(\tsr J∇H)\tp\ \tsr J\ (\tsr J∇H)= \tsr J.
	\label{crl0}
\end{equation}

When the series is truncated at order $N$, the flow of $\breve{u}^N$ is
\[ φ_{t,\breve u^N}:\quad u_0\quad\mapsto\quad \breve u^N(t)=\sum_{n=0}^Nt^nu_n.\]
The following theorem shows that the scheme based on the truncated series is symplectic at order $N+1$.
\begin{theorem}\label{n+1symplectic}
	If the series is convergent then
	\begin{equation}
		(∇φ_{t,\breve u^N})\tp\ \tsr J\ ∇φ_{t,\breve u^N}=\tsr J\ +\ O(t^{N+1})
	\end{equation}
for $t\in[0,\delta t]$ where $\delta t$ is the convergence radius. 
\end{theorem}
Indeed, 
\begin{align*}
	(∇φ_{t,\breve u^N})\tp\ \tsr J\ ∇φ_{t,\breve u^N}&=∑_{n=0}^N∑_{k=0}^n÷{t^n}{k!(n-k)!}(\tsr JD^{k}∇H)\tp\ \tsr J\ (\tsr JD^{n-k}∇H)
	\\&+∑_{n=N+1}^{2N}∑_{k=n-N}^N÷{t^n}{k!(n-k)!}(\tsr JD^{k}∇H)\tp\ \tsr J\ (\tsr JD^{n-k}∇H).
\end{align*}
Using (\ref{crl}) and (\ref{crl0}), the theorem follows. Note that in theorem \ref{n+1symplectic}, $δt$ is generally small.

\vspace{\baselineskip}

In the following subsections, BPL is implemented and tested on a Hamiltonian equation. Next, we present some experiments on non-Hamiltonian equations.

In simulations, the truncation order of the series is set to $N=10$ unless otherwise stated. The degree of the numerator and the denominator of the Padé approximant are $N_{num}=4$ and $N_{den}=5$. A singular value decomposition is used to strengthen the robustness of the Padé calculation \cite{gonnet11}. Twenty Gauss-Laguerre roots are used for the quadrature.

The aim of these simulations is not to make an extensive comparison of BPL with classical schemes (this will be done in a forthcoming paper) but only to show the potential of the scheme in predicting long time dynamics.

\subsection{Periodic Toda lattice}

We consider again the periodic Toda lattice from section \ref{section:toda}. The quality parameter $ε_{res}$ of BPL is choosen such that the mean time step $δt$ is approximately 0.1, and compare the results with that of RK4 and RK4sym (see figures \ref{fig_toda1_errH} and \ref{fig_toda1_eigen}).

\begin{figure}[ht]
	\includegraphics[width=6cm]{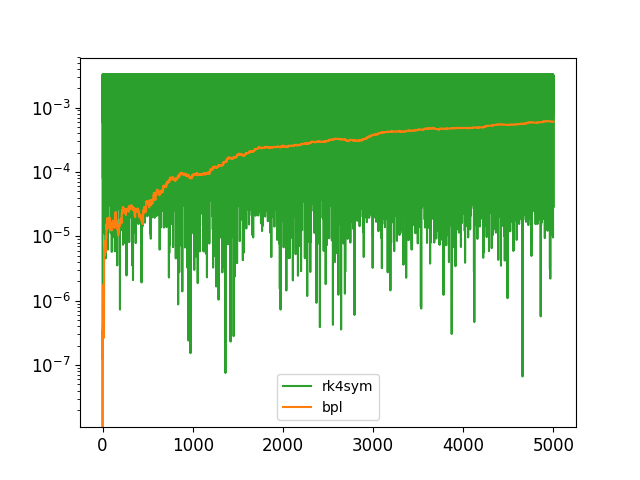}
	\includegraphics[width=6cm]{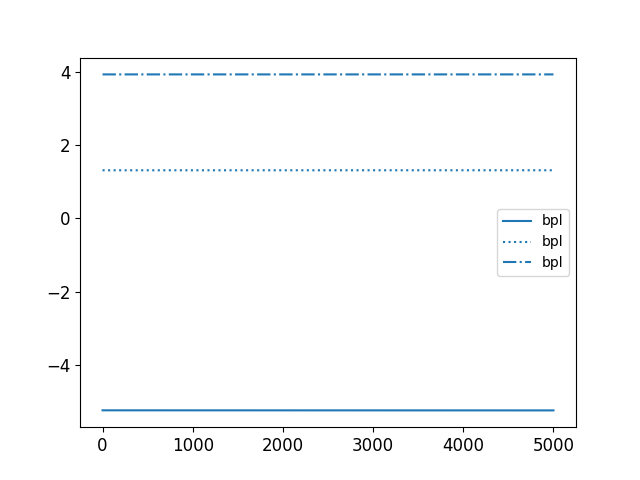}
	\caption{Toda lattice with $δt\simeq 0.0983$. Left: Local errors on $H$. Right: Eigenvalues computed with BPL.}
	\label{fig_toda1_bpl_errH}
\end{figure}

Figure \ref{fig_toda1_bpl_errH}, left, presents the local relative errors on the Hamiltonian. As can be seen, BPL is much more accurate than RK4sym for $t\in[0,5000]$. The value of the global error, defined as
\[ E_H^{mean}= ÷1{t_f}\int_0^{t_f}÷{|H(\q,\p)-H(\q^0,\p^0)|}{|H(\q^0,\p^0)|} ,\quad\quad t_f=5000,\]
is $3.010·10^{-4}$ with BPL and $2.261·10^{-3}$ with RK4sym. BPL also preserves the eigenvalues of the matrix $L$ in equation (\ref{l}), as seen on the right part of figure \ref{fig_toda1_bpl_errH}.

%
%

Table \ref{sametimestep} compares the CPU time needed for 5000 seconds of simulation. It shows that RK4 is the fastest but, as already mentionned, it is not accurate enough for $Δt=0.1$. It also shows that BPL, for approximately the same mean time step, is about twice as slow as RK4sym, but is 7.5 times more accurate.

\begin{table}[ht]
	\begin{tabular}{|c|c|c|c|}
		\hline
		&Mean time-step&CPU&Mean error\\
		\hline
		RK4&0.1&107.44&$3.7902.10^{-1}$
		\\
		RK4sym&0.1&128.74&$2.261·10^{-3}$
		\\
		BPL   &0.0983&259.31&$3.010·10^{-4}$
		\\
		\hline
	\end{tabular}
	\caption{Toda lattice. CPU and accuracy comparison, with almost the same (mean) time step.}
	\label{sametimestep}
\end{table}


In a second test, the different parameters (time step for RK4 and RK4sym and $ε_{res}$ for BPL) are set such that the global accuracies are comparable. Table \ref{sameerror} shows the (mean) time steps and the CPU errors for $E_H^{mean}$ around $2.43·10^{-3}$. As can be seen, RK4sym needs 28 percent less time than BPL to achieve the same accuracy on $H$, due to its especially good property towards the preservation of the Hamiltonian. However, BPL is more than 8 times as fast as the classical RK4 scheme. 

\begin{table}[ht]
	\begin{tabular}{|c|c|c|c|}
		\hline
		&Mean time-step&CPU&Mean error\\
		\hline
		RK4&0.0275&1475.49&$2.130·10^{-3}$
		\\
		RK4sym&0.1&128.74&$2.261·10^{-3}$
		\\
		BPL   &0.125&179.08&$2.897·10^{-3}$
		\\
		\hline
	\end{tabular}
	\caption{Toda lattice. Time step and CPU comparison, with almost the same mean error on $H$. Final time: 5000.}
	\label{sameerror}
\end{table}

\subsection{Duffing equation}

In the next numerical experiment, consider the forced Duffing equation
\begin{equation}
	\ddot u + r\dot u+au+bu^3=c\cos(ωt)
	\label{duffing}
\end{equation}
which describes nonlinear damped oscillators \cite{thompson02,jordan07}. 
We first consider the force-free case with two sets of coefficients for which there is a first integral:
\begin{itemize}
	\item {\it Case 1}:  $a=2/9,\,b=1,\,r=-1$,
	\item {\it Case 2}:  $a=1,\,b=1,\,r=0$.
\end{itemize}
In {\it Case 1}, the first integral is
\begin{equation}
	I=\e^{-÷{4t}3}\left(\dot u^2-÷23u\dot u +÷19u^2+÷12u^4\right).
	\label{first1}
\end{equation}
In {\it Case 2}, the equation can be written in a Hamiltonian form, with a Jacobi elliptic sine function as exact solution. The first integral (the Hamiltonian function) is
\begin{equation}
	÷12\dot u^2+÷12 u^2+÷14u^4.
	\label{first2}
\end{equation}
The initial conditions are $u(0)=1$, $\dot u(0)=0$.

\begin{figure}[ht]
	\includegraphics[width=6cm]{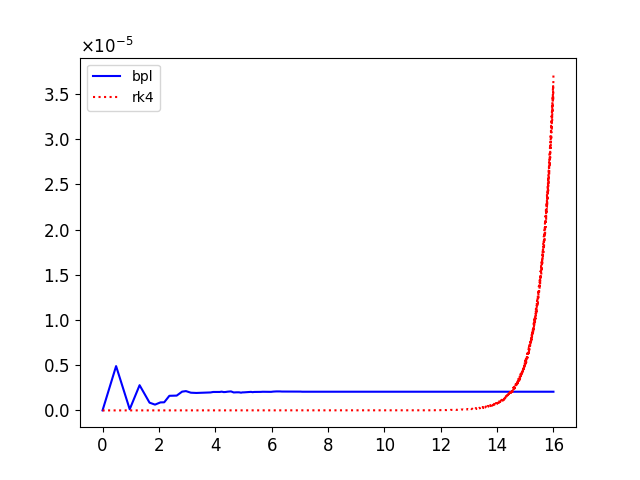}
	\caption{Duffing equation, {\it Case 1}. Relative error on the first integral  (\ref{first1}).}
	\label{fig_duffing_case1_err}
\end{figure}

In {\it Case 1}, the quality criterium $ε_{res}$ of BPL is set to $10^{-4}$. The mean time step within the first 20 seconds is approximately $4.7·10^{-4}$. 
Figure \ref{fig_duffing_case1_err} plots the evolution of the relative error on the first integral, compared to the relative error of RK4 with a time step $4·10^{-4}$. As can be seen, the RK4 error is very small in the beginning of the simulation. It remains below $10^{-9}$ at $t=10$. But when $t$ is large, the RK4 error becomes big and reaches 50 percent at $t=22.3s$. With BPL, the error oscillates when $t$ is small. The amplitude is of the same order as $ε_{res}$. Next, the error stabilizes rapidly around $1.27·10^{-6}$.

\begin{figure}[ht]
	\includegraphics[width=6cm]{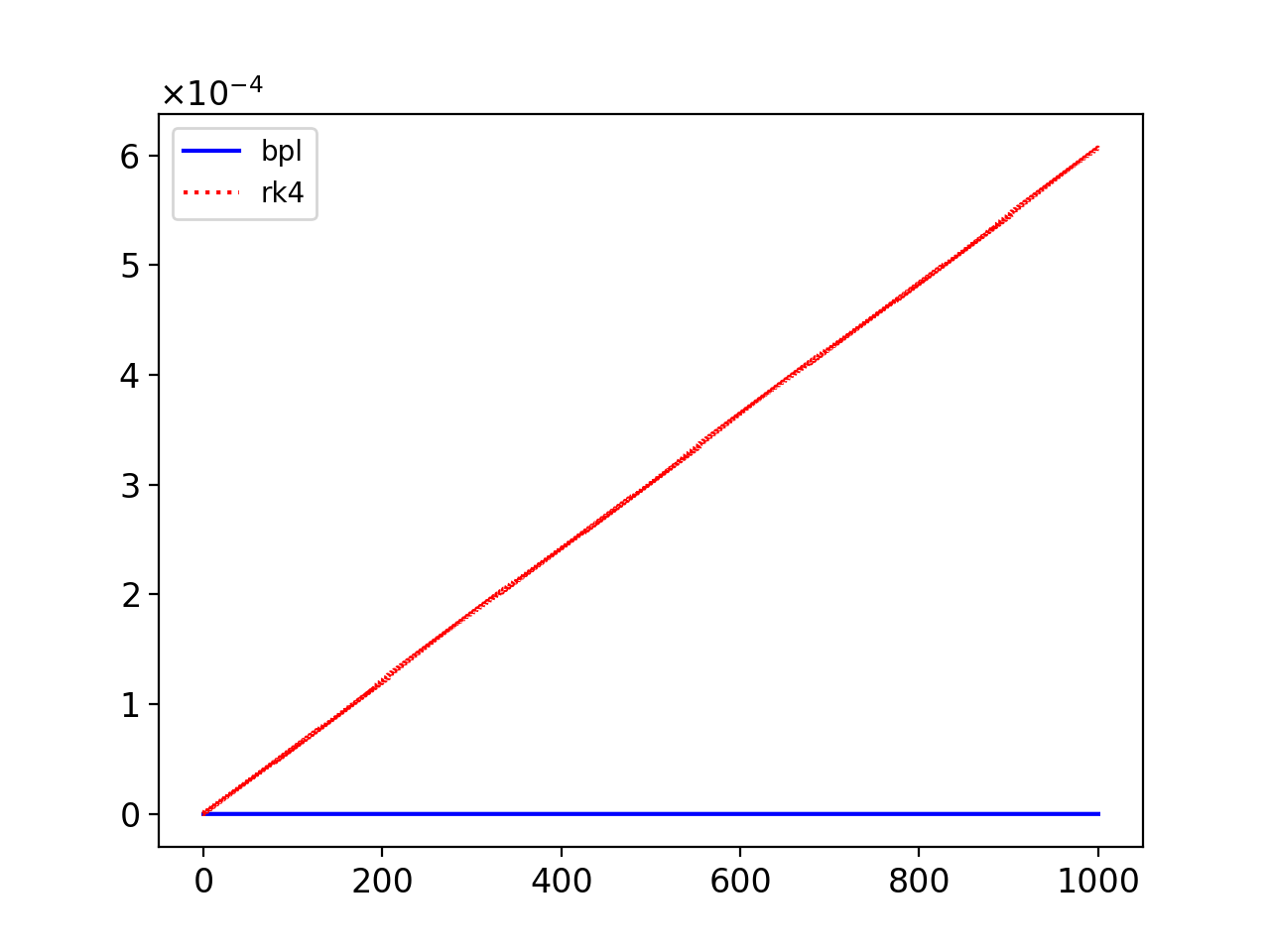}
	\caption{Duffing equation, {\it Case 2}. Relative error on the first integral  (\ref{first2}).}
	\label{fig_duffing_case2_err}
\end{figure}

In {\it Case 2}, the quality criterium $ε_{res}$ of BPL is choosen such that the mean time step is around $0.106$. For RK4, the time step is set to 0.1. The evolution of the relative error on the first integral (\ref{first2}) is plotted on figure \ref{fig_duffing_case2_err}. As can be seen, the error of BPL is much smaller than that of RK4.

\begin{figure}[tp]
	\includegraphics[width=6cm]{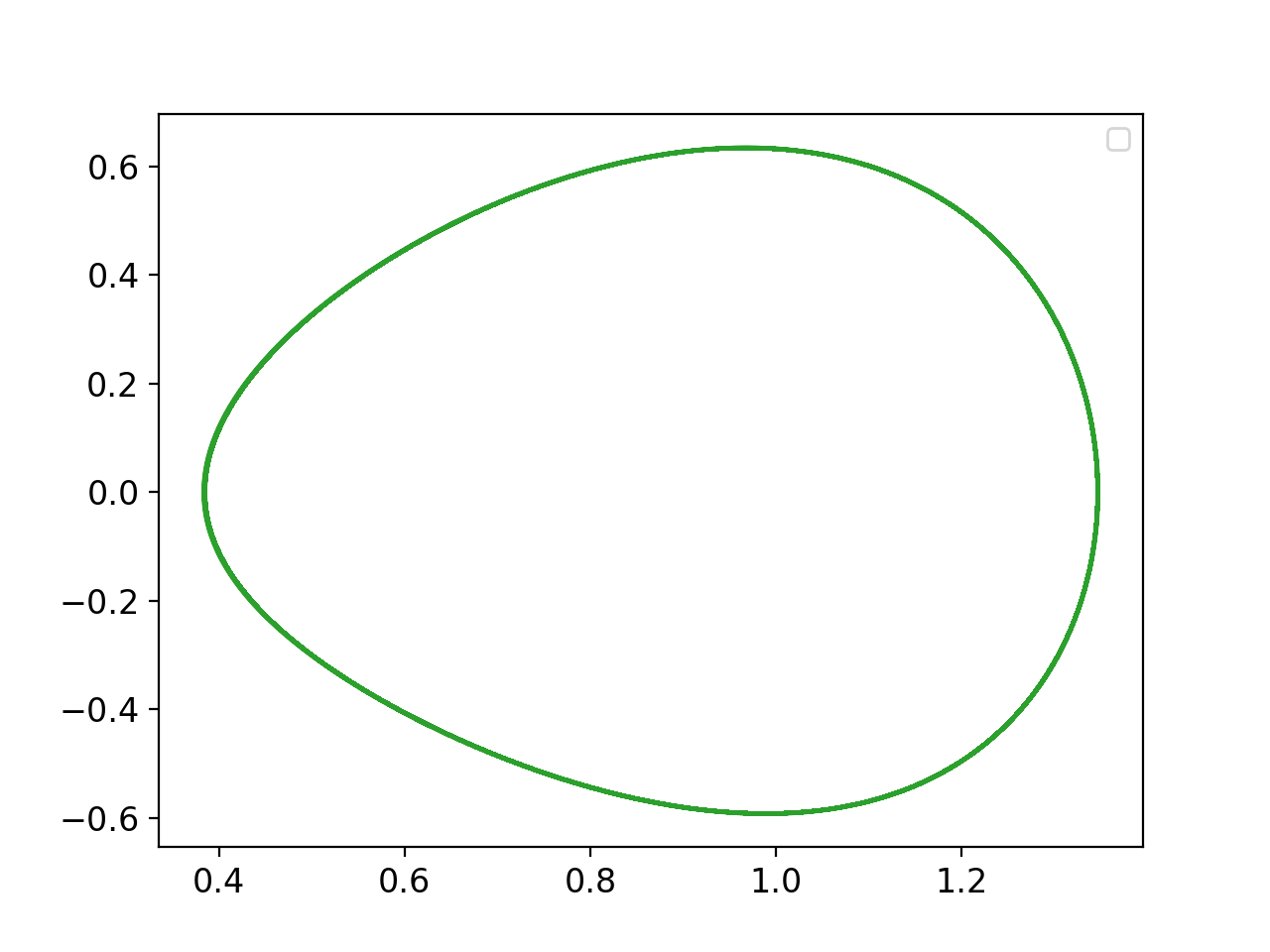}
	\includegraphics[width=6cm]{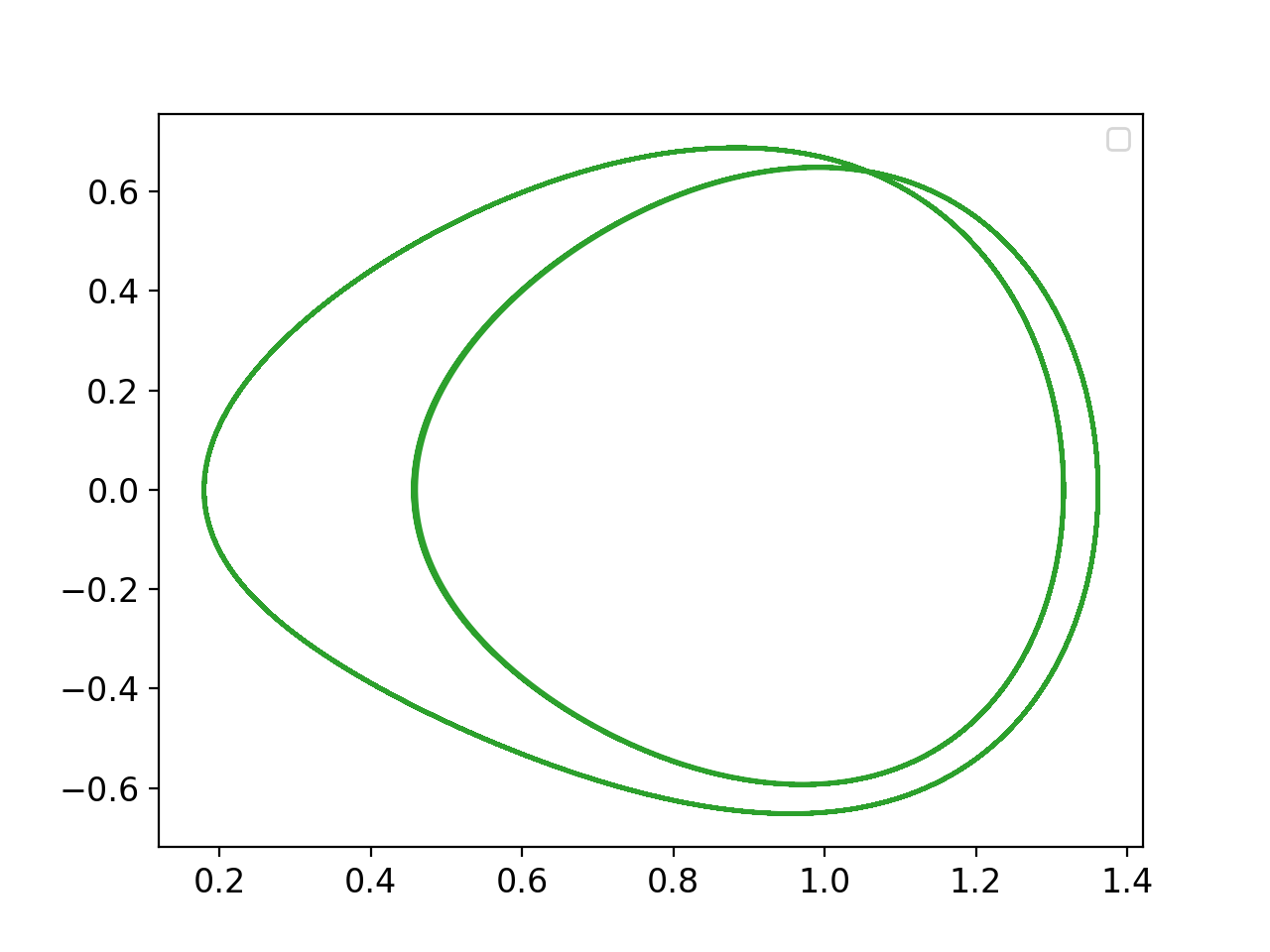}
	\includegraphics[width=6cm]{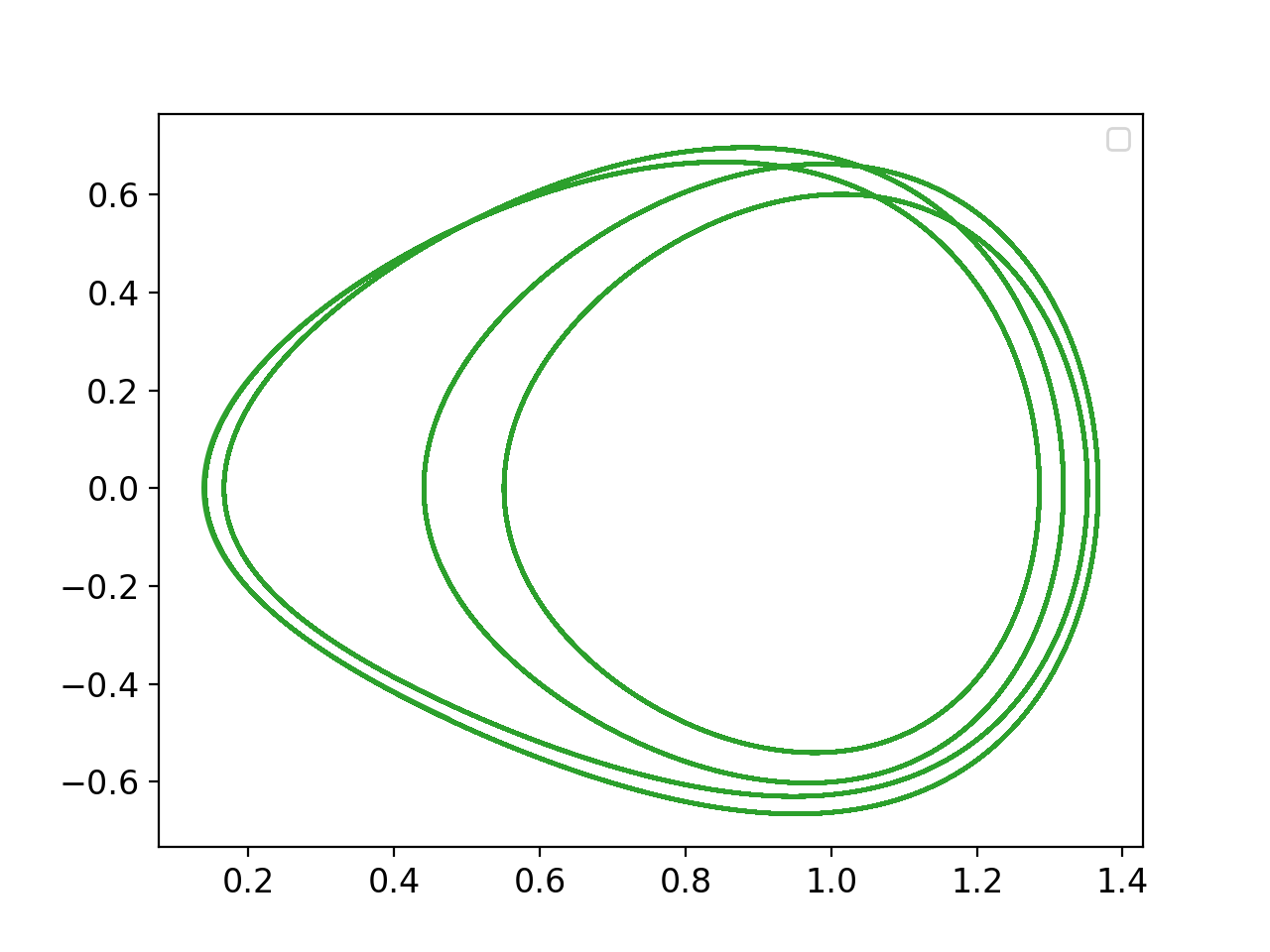}
	\includegraphics[width=6cm]{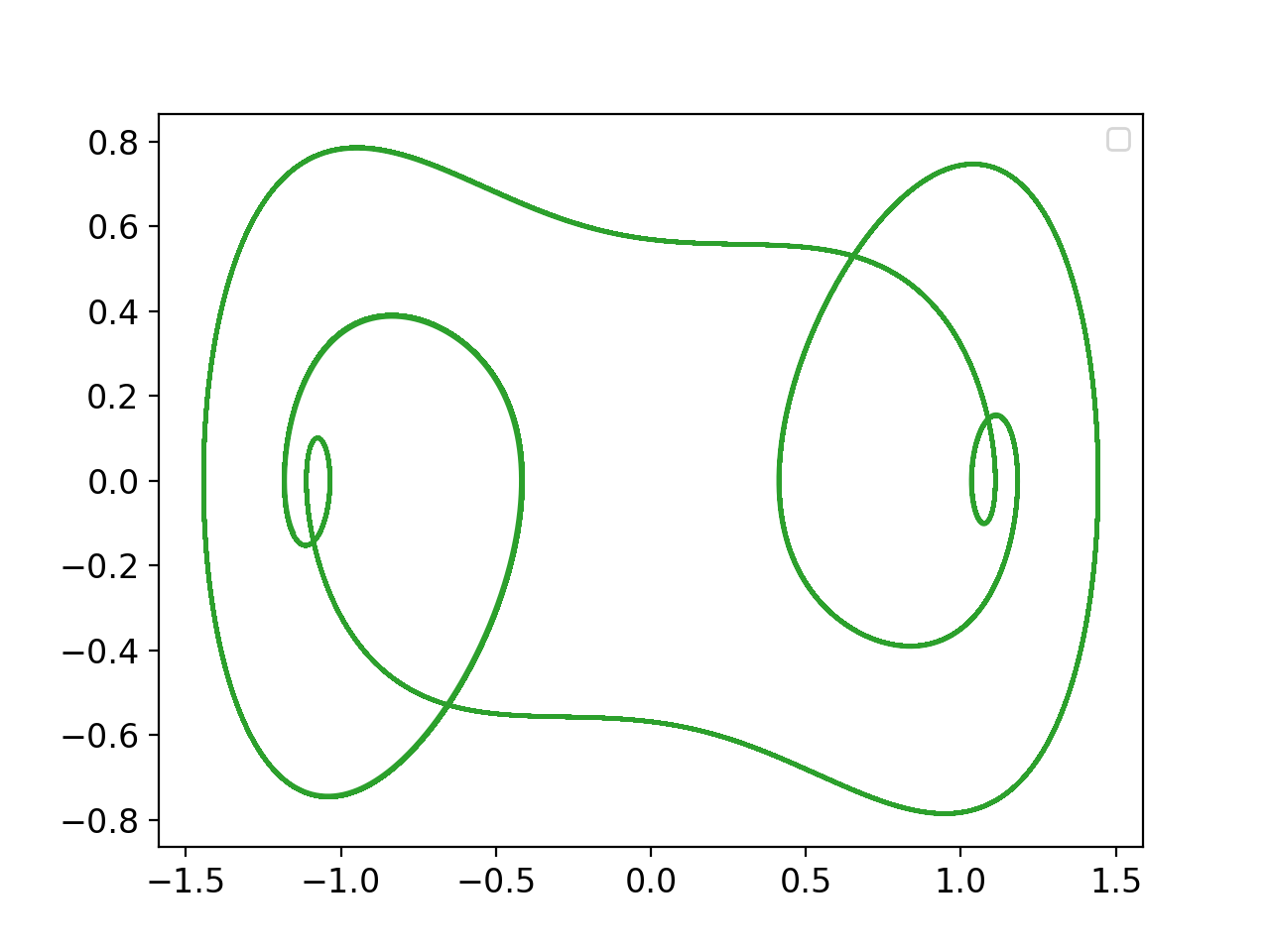}
	\includegraphics[width=6cm]{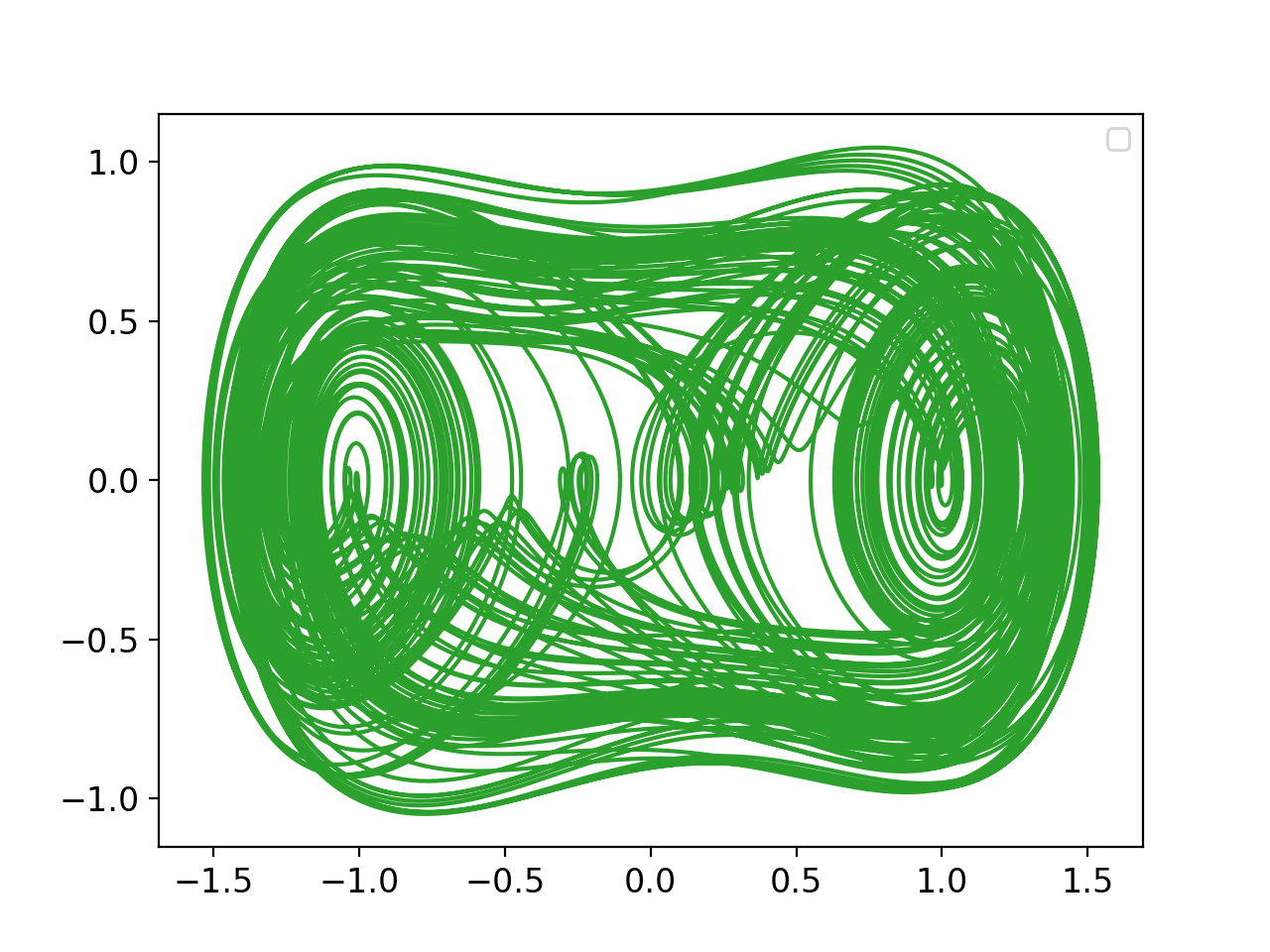}
	\includegraphics[width=6cm]{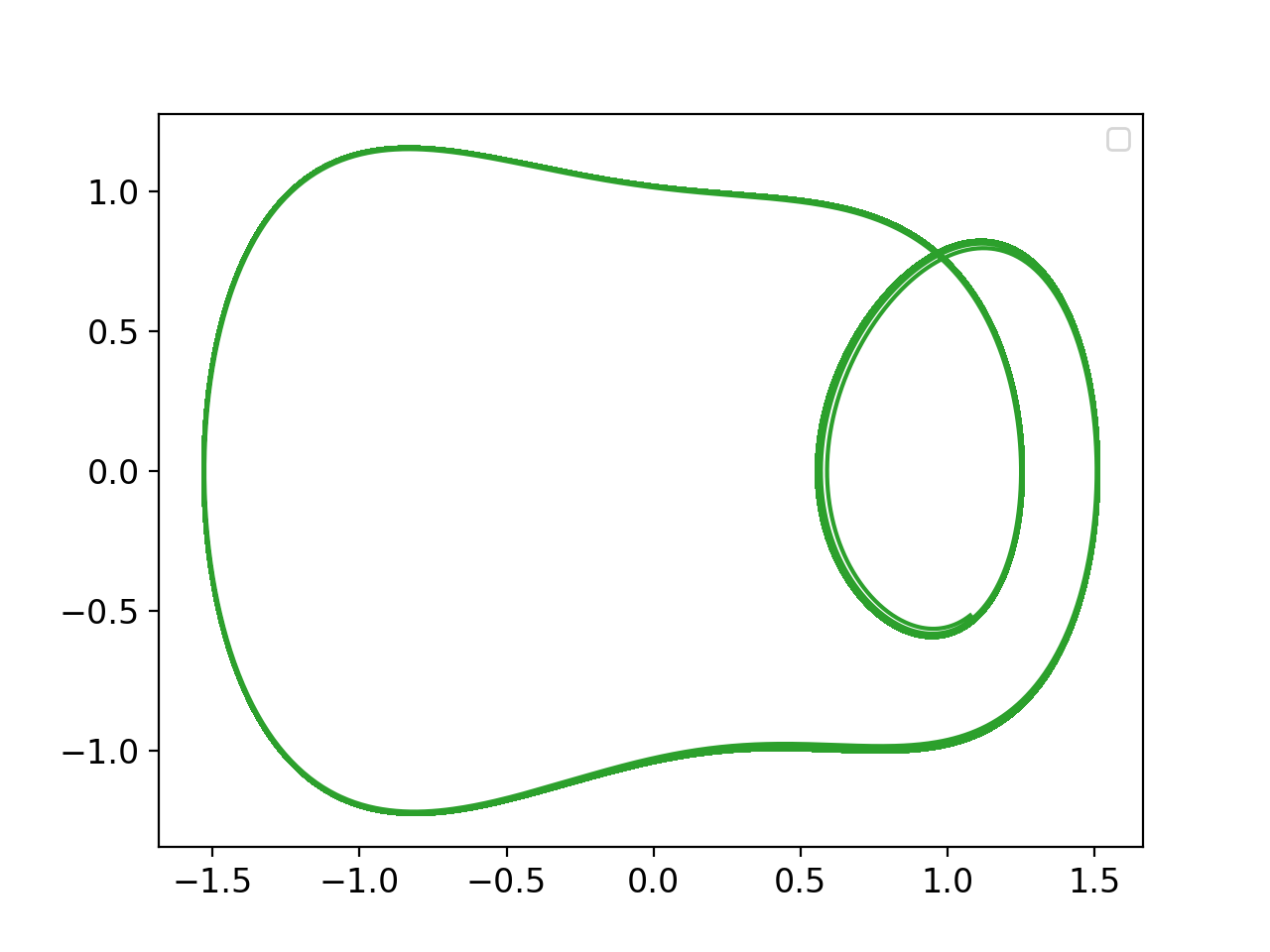}
	\caption{Duffing equation: phase portraits. From left to right and from top to bottom, $c=0.20$, $0.28$, $0.29$, $0.37$, $0.50$, $0.65$ }
	\label{fig_duffing_c}
\end{figure}

To end up with Duffing equation, some phase portraits obtained with BPL are computed. They corresponds to $a=-1$, $b=1$, $r=.3$, $ω=1.2$ and $c$ varying from 0.20 to 0.65. The initial conditions are $u(0)=1$ and $\dot u(0)=0$. Figure \ref{fig_duffing_c} presents the phase trajectories for $t\in[40,1000]$, that is after the transient phase. These plots have been obtained with a rather loose value of $ε_{res}$ for which the mean time step is around 0.5. But as can be seen, when $c=0.20$, $0.28$, $0.29$, $0.37$, and $0.65$, the (multiple)-periodicity is very well captured even over a very long time interval. Indeed, the curves are closed. For $c=0.50$, the solution is chaotic but is bounded. These results are in agreement with those presented in \cite{jordan07}.

\vspace{\baselineskip}

In the last subsection, BPL is applied to a semi-discretized partial differential equation. It is compared to some other adaptative schemes. Since the system is big enough, it is worth to give an indication on the CPU simulation time.

\subsection{Korteveg-de-Vries equation}

Consider the Korteweg-de-Vries equation 
\begin{equation}
\label{kdv}
\cfrac{\partial u}{\partial t}  + c_0 \cfrac{\partial u}{\partial x} + \beta \cfrac{\partial^3 u}{\partial x^3} + \cfrac{\alpha}2 \cfrac{\partial u^2}{\partial x}= 0
\end{equation}
which models waves on shallow water surfaces \cite{korteweg95}. In this equation, the linear propagation velocity $c_0$, the non-linear coefficient $\alpha$ and the dispersion coefficient $\beta$ are positive constants, linked to the gravity acceleration $g$ and the mean depth $δ$ of the water by
$c_0=\sqrt{gδ}$, $\alpha=\frac32\sqrt{g/δ}$ and $\beta=d^2c_0/6$.

The solution is assumed to be periodic with period $X$ in space. Equation (\ref{kdv}) is then discretized in space with a spectral method. The solution is approximated by its truncated Fourier series:
\begin{equation}
	u(x,t)\simeq \sum_{|m|\leq M}\hat u^m(t)\e^{im\omega x},
	\label{fourier}
\end{equation}
where $M\in \mathbb N$ and $\omega=\frac{2\pi}X$.

With BPL, the Fourier coefficient array $\hat u(t)$ is decomposed into a time series
\begin{equation}
	\hat u^m(t)=\sum_{k=0}^K\hat u^m_kt^k
	\label{bpl_fourier}
\end{equation}
The series coefficients are computed explicitely as follows:
\begin{equation}
	\hat u_{k+1}=\cfrac1{k+1}\left[(-c_0i\omega m+i\beta\omega^3 m^3)\hat u_k+\cfrac12\,i\alpha m\omega \ \sum_{n=0}^k\hat u_n*\hat u_{k-n}\right].
	\label{uk1}
\end{equation}
The initial condition is the periodic prolongation of the function
\begin{equation}
	u_0(x)=h\operatorname{sech}^2(\kappa x), \quad\quad\quad x\in\left[ -\frac X2,\frac X2 \right]
	\label{kdv_initial}
\end{equation}
with 
$\kappa =\sqrt{3h/4δ^3}$.
The exact solution is 
\begin{equation}
	u(x,t)=u_0(x-ct)
	\label{kdv_solution}
\end{equation}
where $c=c_0(1+h/2δ).$
We take $X=24\pi$, $δ=2$, $g=10$ and $h=\frac12$. The period is $T\simeq14.98$s. To begin with, the size of the system is set to $d=128$ (that is the number of spectral discretization points is 129).

BPL is compared to two other schemes. The first one is the adaptative 4-th order Runge-Kutta scheme (still denoted RK4 in this subsection). This scheme is explicit. The second one is the exponential time differencing associated to RK4 (denoted ETDRK4), developed by Cox and Matthews in \cite{cox_2002}. This scheme is based on an exact, exponential type, resolution of the linear part of the equation, followed by an explicit adaptative Runge-Kutta resolution of the non-linear part. The algorithm is not completely explicit since it requires the (pseudo-)inversion of a matrix. Moreover, it generally needs the evaluation of a matrix exponential, which is numerically expensive. This evaluation is done via a Padé approximants in simulations.

\begin{figure}[htp]
	\centering
	\includegraphics[width=.6\textwidth]{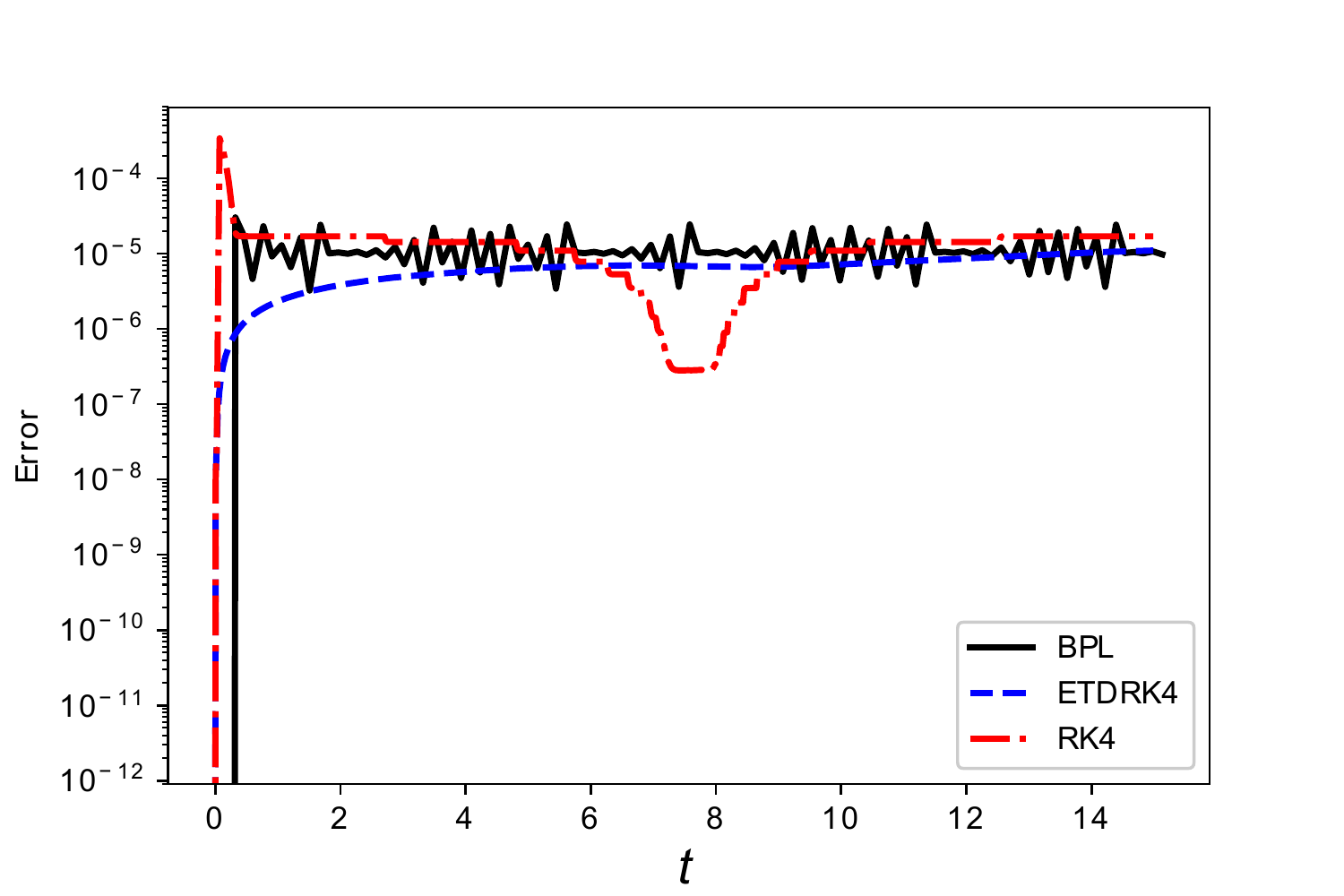}
	\caption{Korteweg de Vries equation. Evolution of the error with time}
	\label{fig:kdv_error_128}
\end{figure}

The precision criteria are calibrated such that the \textit{a posteriori} errors of the three schemes have approximately the same magnitude, as can be seen in figure \ref{fig:kdv_error_128}. This figure shows the difference between the predicted solutions with (\ref{kdv_solution}).

\begin{figure}[htp]
	\centering
	\includegraphics[width=.6\textwidth]{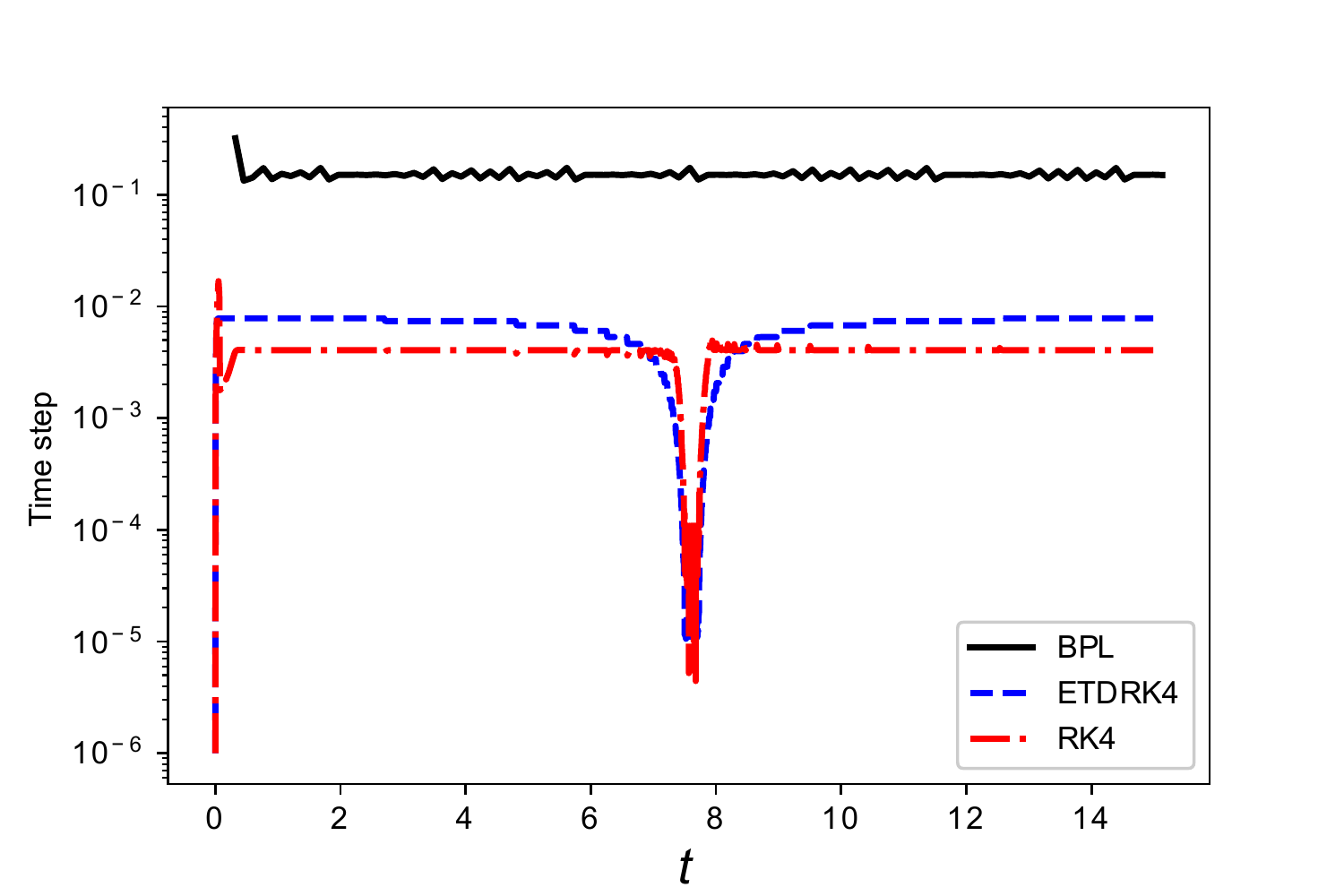}
	\caption{Evolution of the time step with time}
	\label{fig:kdv_timestep_128}
\end{figure}

The time steps are presented in figure \ref{fig:kdv_timestep_128}. In mean, the BPL time step is 238 times bigger than that of ETDRK4. This reflects on the CPU time. Indeed, as can be observed on table \ref{kdv_128}, BPL is about 950 times faster than ETDRK4 for approximately the same precision.
Note that, for this specific problem, the time step with RK4 has the same order as that of ETDRK4, but RK4 is more efficient than ETDRK4 in terms of computation time, since it requires neither numerical matrix (pseudo-)inversion nor exponential. Compared to BPL, RK4 takes 60 times more CPU time to reach one period.

\begin{table}[ht]\centering
\begin{tabular}{|c|c|c|c|}
\hline
&  BPL & ETDRK4 & RK4 \\ 
\hline
Mean time step & 1.53$\,\cdot 10^{-01}$ &6.41$\,\cdot 10^{-04}$ &1.88$\,\cdot 10^{-03}$  \\
\hline
$L^2$ error at $t=T$& 9.76$\,\cdot 10^{-06}$ & 1.10$\,\cdot 10^{-05}$ &1.69$\,\cdot 10^{-05}$  \\
\hline
Simulation time & 1.74 &1.66$\,\cdot 10^{+03}$  &5.98$\,\cdot 10^{+01}$ \\
\hline
\end{tabular}
\caption{Time step, error and CPU time over one period, with $d=128$}
\label{kdv_128}
\end{table}

\begin{figure}[htp]
	\centering
	\includegraphics[width=.6\textwidth]{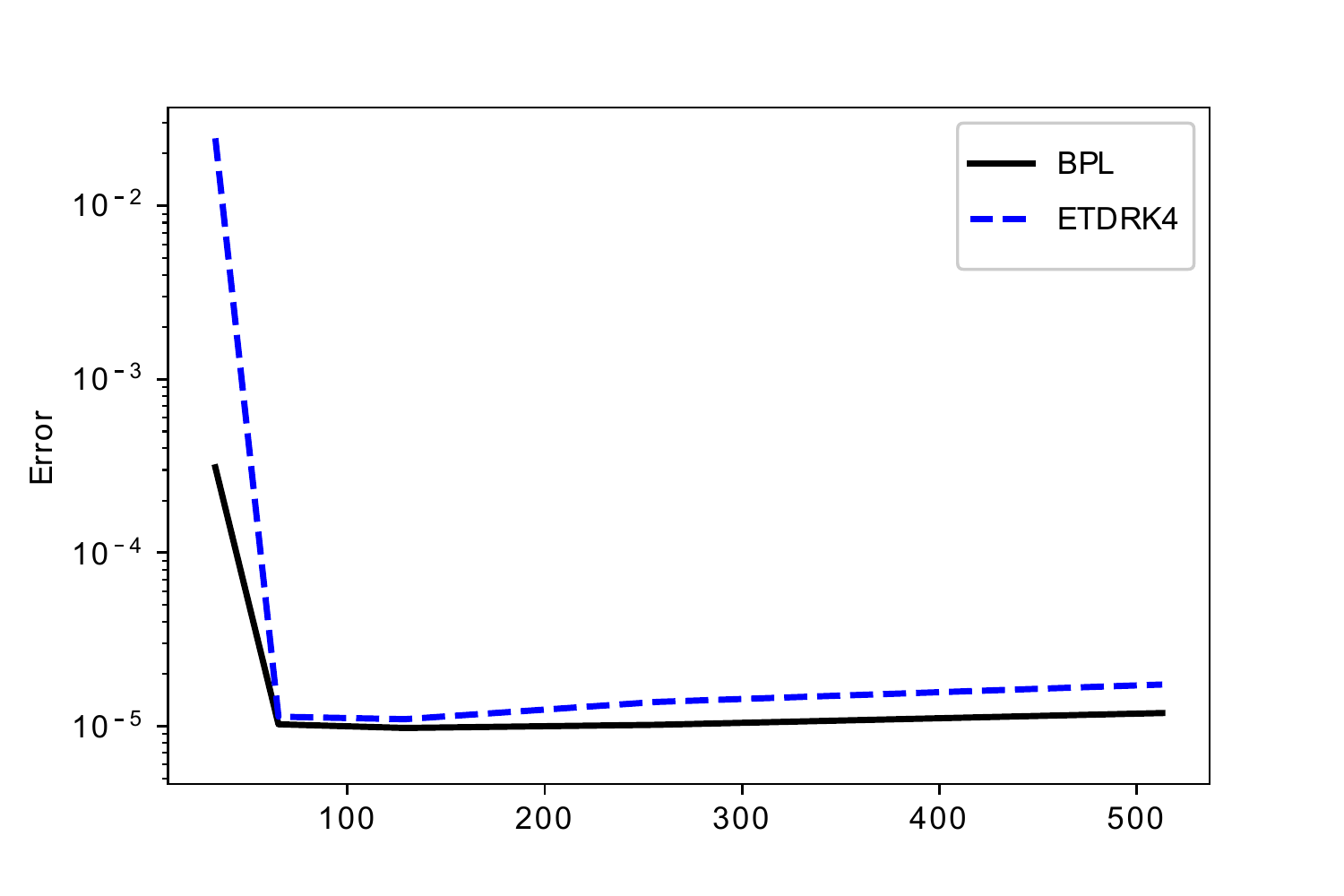}
	\\[-.5\baselineskip]
	a) $L^2$ error when $t=T$ 

	\includegraphics[width=.6\textwidth]{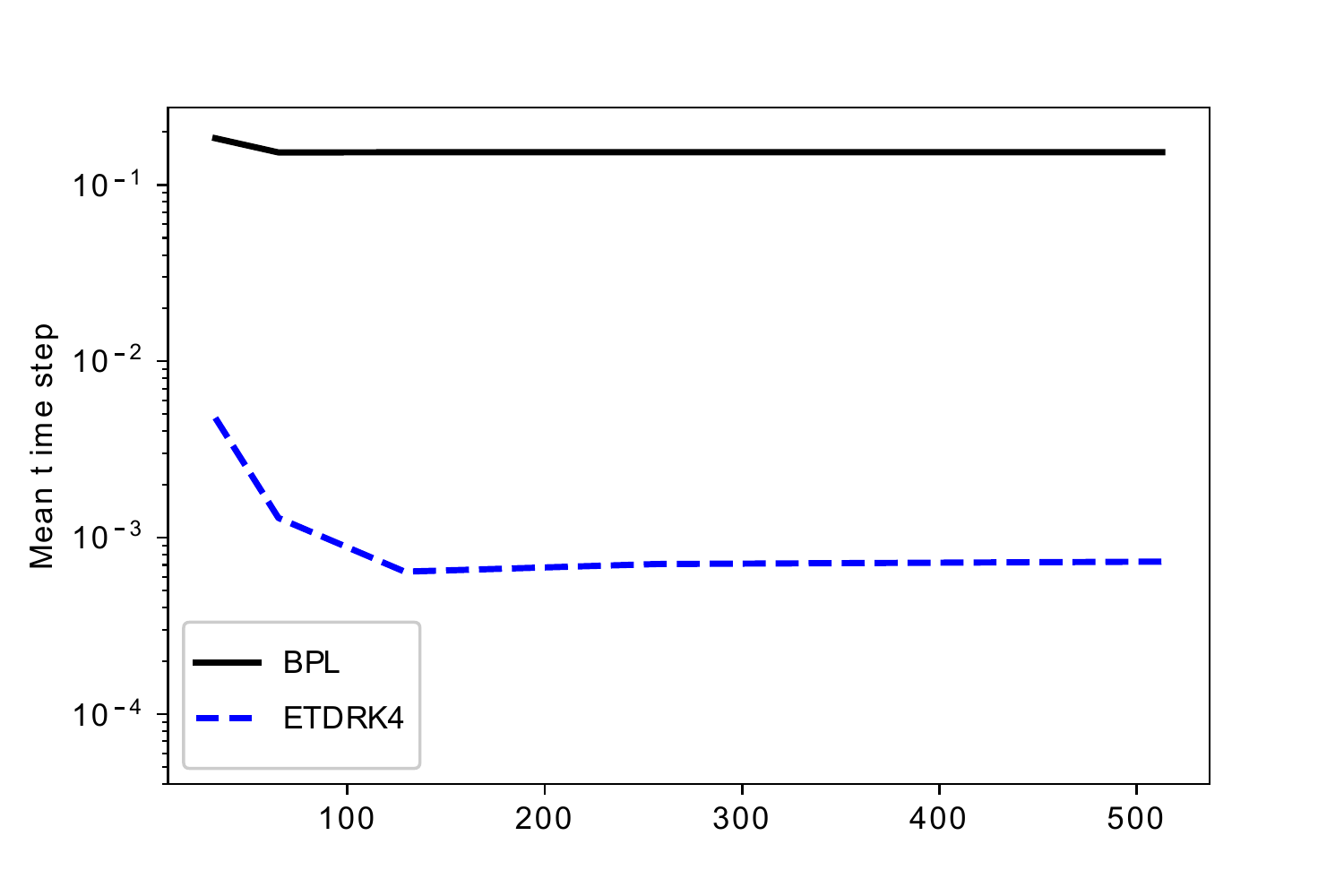}
	\\[-.5\baselineskip]
	b) Mean time step over one period
	
	\includegraphics[width=.6\textwidth]{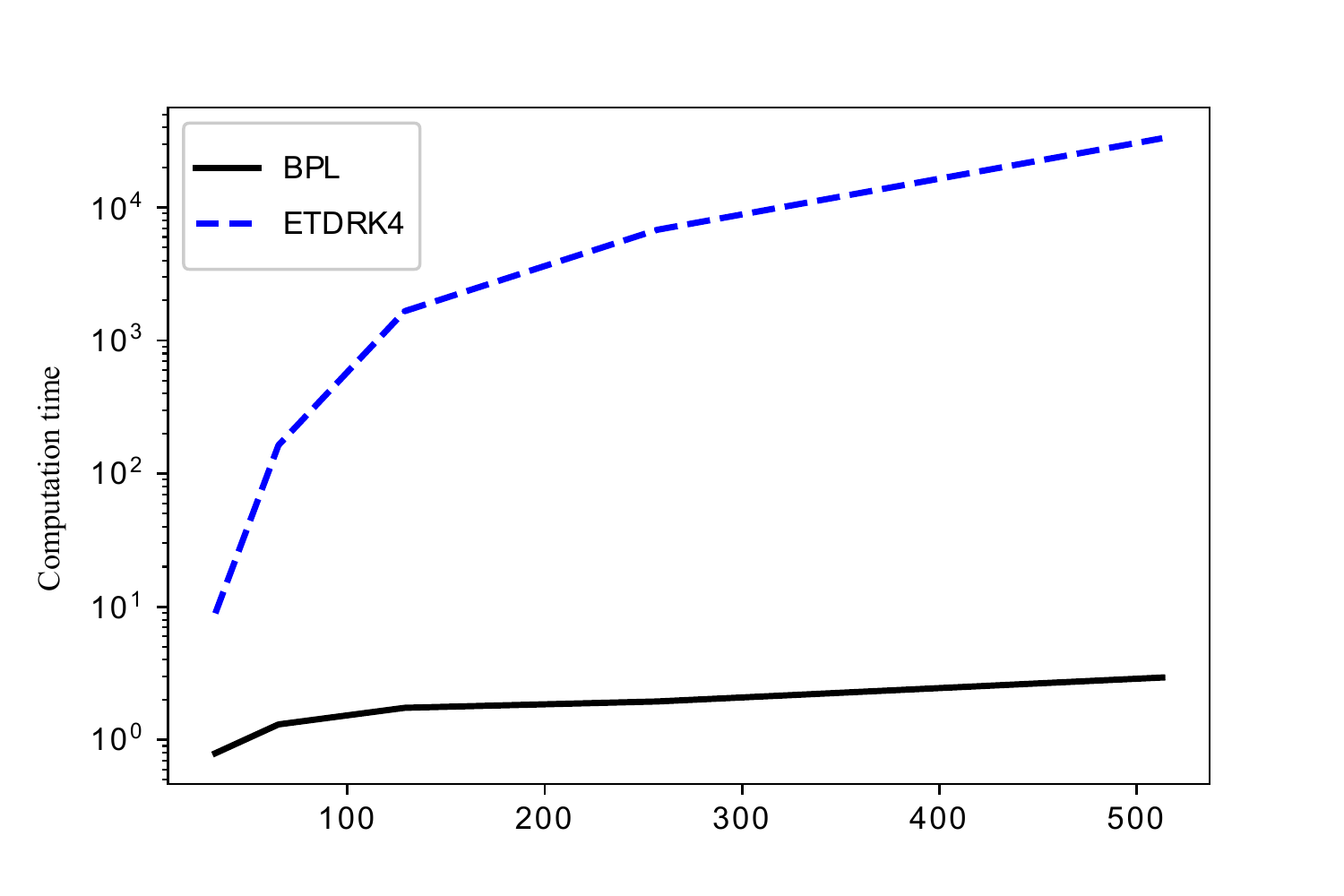}
	\\[-.5\baselineskip]
	c) CPU time over one period
	\caption{Evolution of the error, the time step and the computational time with the size $d$ of the problem}
	\label{fig:kdv}
\end{figure}

In the next simulation, we analyse the behaviour of the schemes when the size $d$ of the problem is increased. Figure \ref{fig:kdv}a presents the $L^2$ error for $t=T$. It shows that the precision of BPL and ETDRK4 remains approximately the same, except when $d$ is very small. Figure \ref{fig:kdv}b shows however that BPL requires much less iterations (100 iterations versus 20389 for ETDRK4 when $d=512$ to reach one period). 
The time steps of BPL and ETDRK4 seem to have the same behaviour when $d$ is large enough. Indeed, they tend to be independent of $d$, as suggested by figure \ref{fig:kdv}b. However, the computation time increases much more rapidly with ETDRK4. For BPL, the growth rate of the CPU time between $d=128$ and $d=512$ is 51 percent whereas, for ETDRK4, it is 381 percent.

\begin{figure}[htp]
	\centering
	\includegraphics[width=.6\textwidth]{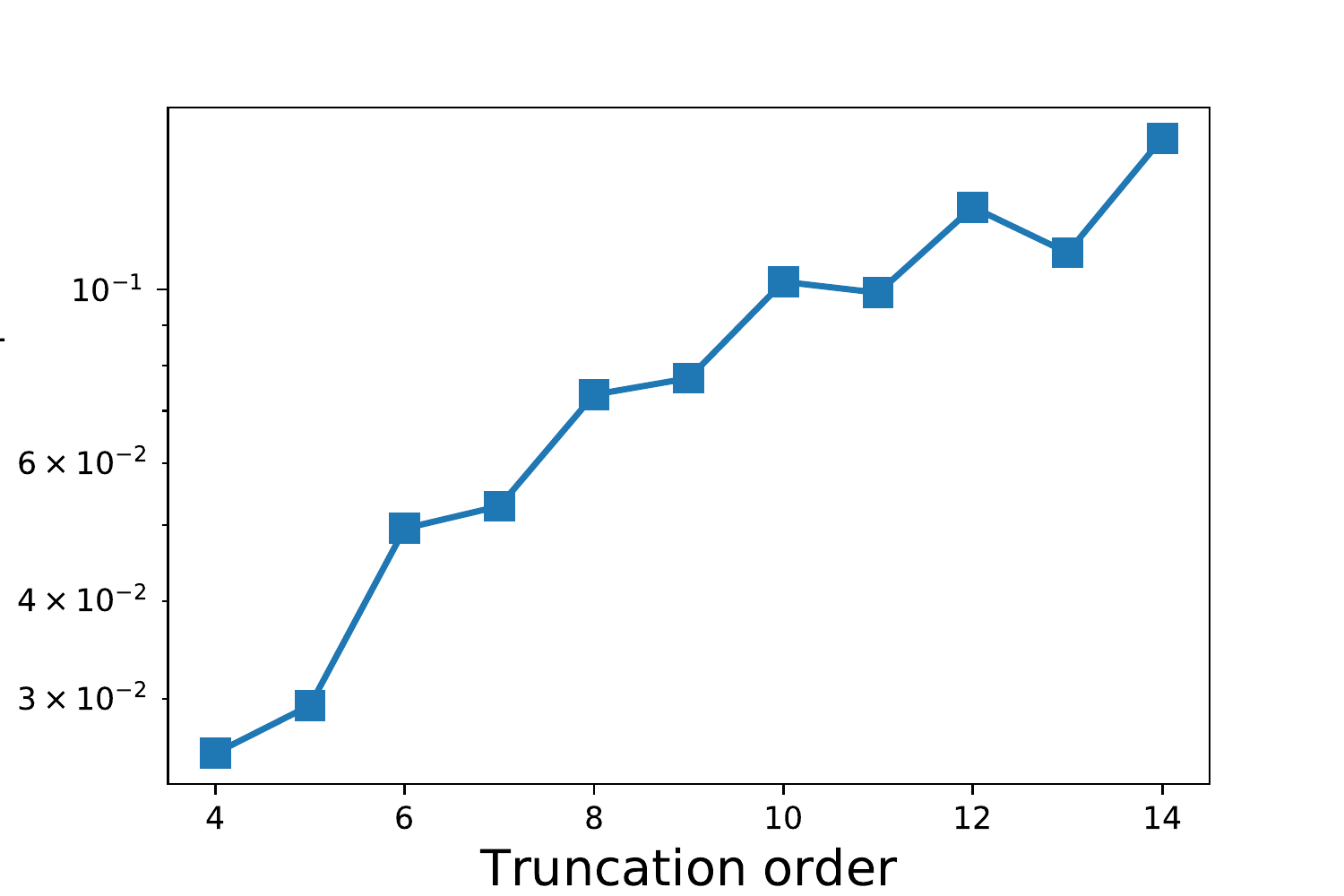}
	\\
	a) Mean time step

	\includegraphics[width=.6\textwidth]{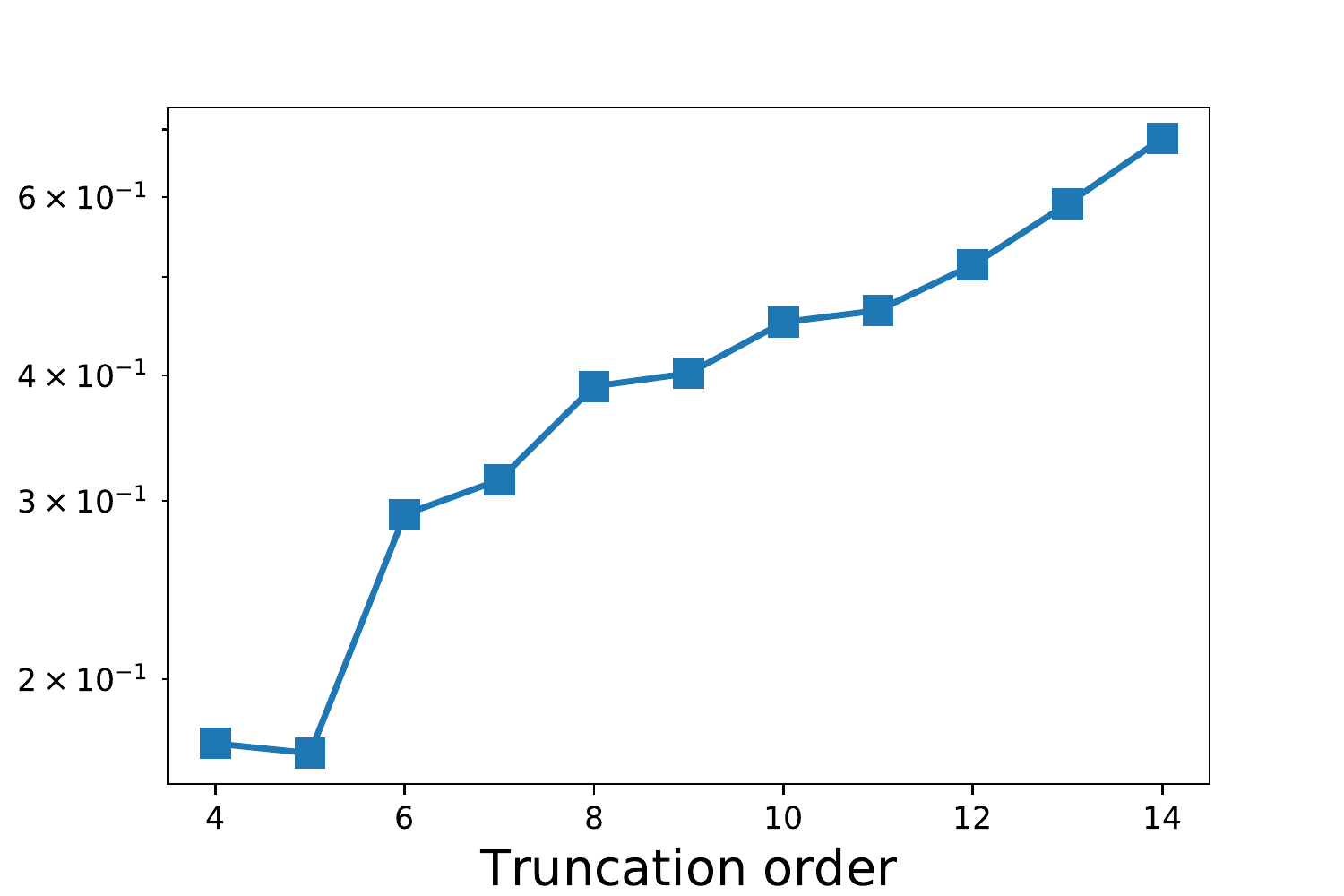}
	\\
	b) CPU time

	\includegraphics[width=.6\textwidth]{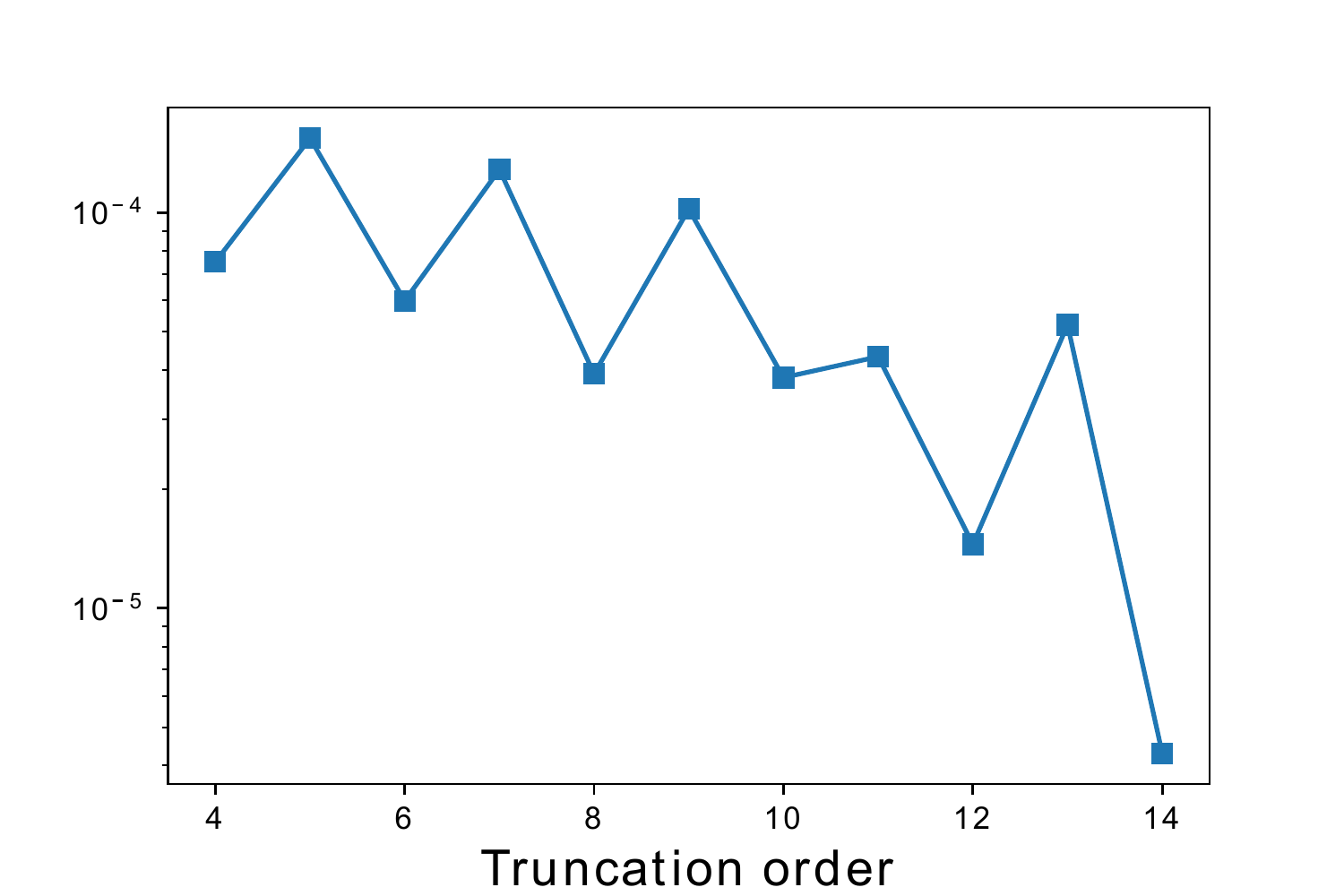}
	\\
	c) $L^2$ error when $t=T$
	\caption{Evolution of the time step, computation time and error with the order $N$ of truncature of the time series in BPL}
	\label{fig:kdv_order}
\end{figure}

In all of the previous simulations, the truncation order $N$ of the time series in BPL was set to 10. In our last test, the effect of $N$ on the quality of BPL is analysed. For this, the size of the problem is fixed to $d=128$. Figure \ref{fig:kdv_order}a shows that the time step increases with $N$, passing from $Δt_{mean}=0.0256$ for $N=4$ to $Δt_{mean}=0.156$ when $N=14$. Despite the number of iterations is consequently reduced, the CPU time also increases with $N$, going from 0.686$s$ to 0.173$s$, as can be observed in figure \ref{fig:kdv_order}b. This is caused by the fact that more coefficients of the series and more Padé coefficients have to be computed. As for it, the error fluctuates but globally decreases from $7.51\cdot 10^{-5}$ to $4.30\cdot 10^{-6}$. This fluctuation is not uncommon in series based approximations. It is interesting to note that whereas the error is divided by $17.5$, the CPU time is multiplied only by $3.96$ between $N=4$ and $N=14$. In other words, the precision increases faster than the CPU time when the order of the scheme is increased.

\section{Conclusion}

In this article, we gave an overview of some time integrators for long-time simulations. Two geometric integrators and a general-purpose time integrator was presented.

Through numerical examples, the ability of symplectic integrator in preserving the Hamiltonian, the angular momentum or eigenvalues was observed. Moreover, it was shown that symplectic integrators are more robust compared to classical schemes when the time step is enlarged (in the example of Toda lattice) or when a perturbation is introduced (three-body problem). 

Next, a way of constructing Dirac integrators for constrained system was given. Numerical experiments showed that respecting the Dirac structure at discrete level avoids numerical artifacts. As a consequence, Dirac integrators are able to reproduce the dynamics of the system over a long time.

Finally, we showed that BPL competes with symplectic integrators in predicting Hamiltonian dynamics (Toda lattice and {\it Case 2} of Duffing equation). For more general equations, BPL also preserves with high precision the first integral of the system, as well as the periodicity when the solution is periodic. Lastly, compared to two popular schemes, BPL appears to require less computation time.

\end{document}